\documentclass[10pt]{iopart}

\usepackage{latexsym}
\usepackage{amsfonts}
\usepackage{epsfig}

\newtheorem{theorem}{Theorem}[section] 
\newtheorem{lemma}{Lemma}[section] 
 
\newtheorem{definition}{Definition}[section] 
\newtheorem{remark}{Remark}[section]

\newcommand{\bfv}{{\bf v}}

\newcommand{\bfu}{{\bf u}}
\newcommand{\bfz}{{\bf z}}
\newcommand{\bfE}{{\bf E}} 
\newcommand{\bfH}{{\bf H}}
\newcommand{\bfx}{{\bf x}} 
 
\newcommand{\bfq}{{\bf q}}
\newcommand{\bfqh}{\hat{{\bf q}}}
\newcommand{\bfd}{{\bf d}} 
\newcommand{\bfn}{{\bf n}} 
\newcommand{\bfr}{{\bf r}} 
\newcommand{\bfs}{{\bf s}} 

\newcommand{\anzus}{J} 
\newcommand{\anzant}{K} 
 
\newcommand{\gammahat}{\hat{\gamma}} 
\newcommand{\Qhat}{\hat{Q}} 
\newcommand{\Yhat}{\hat{Y}} 
\newcommand{\Zhat}{\hat{Z}} 
\newcommand{\Vhat}{\hat{V}} 
\newcommand{\Mhat}{\hat{M}}

\newcommand{\nmax}{\hat{n}} 

\newcommand{\mcalE}{\mathcal{E}}
\newcommand{\hcalE}{\hat{\mathcal{E}}}
\newcommand{\mcalF}{\mathcal{F}}
\newcommand{\mcalJ}{\mathcal{J}}
\newcommand{\mcalT}{\mathcal{T}}
\newcommand{\mcalS}{\mathcal{S}}
\newcommand{\mcalThat}{\hat{\mathcal{T}}}
\newcommand{\mcalU}{\mathcal{U}}
\newcommand{\mcalA}{\mathcal{A}}

\newcommand{\Astar}{A^{*}}
\newcommand{\stilde}{\tilde{\bfs}}
\newcommand{\rtilde}{\tilde{\bfr}}
\newcommand{\scirc}{\hat{\bfs}}
\newcommand{\rcirc}{\hat{\bfr}}

\newcommand{\ndiag}{\mbox{{\rm{\bf diag}}}}
\newcommand{\ngrad}{\mbox{{\rm{\bf grad}}}}
\newcommand{\ncurl}{\mbox{{\rm{\bf curl}}}}
\newcommand{\ndiv}{\mbox{{\rm div}}}

\def\C{{\rm\kern.24em 
   \vrule width.02em 
       height1.4ex depth-.05ex 
   \kern-.26em C}} 
\def\R{{\rm I\kern-.25em R}}

\begin{document}

\title{Time-reversal and the adjoint method with 
an application in telecommunication}

\author{Oliver Dorn}
\address{Departamento de Matem\'{a}ticas, Universidad Carlos III de Madrid, 
Avda. de la Universidad, 30, 28911 Legan\'{e}s, MADRID-SPAIN.} 
\address{Current address: Laboratoire des Signaux et Syst\`{e}mes (L2S),
\'{E}cole Sup\'{e}rieure d'Electricit\'{e} - SUP\'{E}LEC,
Plateau de Moulon,
3, rue Joliot-Curie,
91192 Gif-sur-Yvette Cedex, France.
} 

\ead{oliver.dorn@uc3m.es}

\begin{abstract} 
We establish a direct link between the time-reversal 
technique and the so-called adjoint method for imaging. 
Using this relationship, we derive new solution 
strategies for an inverse problem which arises in telecommunication.
These strategies are all based on iterative time-reversal
experiments, which try to solve the inverse problem
{\em experimentally} instead of computationally. 
We will focus in particular on 
 examples from underwater acoustic communication
and wireless communication in a Multiple-Input Multiple-Output
(MIMO) setup.
\end{abstract}


\maketitle

\section{Introduction}\label{Sec1}

Time-reversal techniques have attracted great attention
recently due to the large variety of interesting potential
applications. The basic idea
of time-reversal (often also referred to as
'phase-conjugation' in the frequency domain) 
can be roughly described as follows:  
A localized source emits a short
pulse of acoustic, electromagnetic or elastic energy which
propagates through a richly scattering environment. 
A receiver array records the arriving waves, typically as 
a long and complicated signal due to the complexity of 
the environment, and stores
the time-series of  measured signals in memory. After a short 
while, the receiver array sends the recorded signals time-reversed
(i.e. first in - last out) 
back into the same medium. Due to the time-reversibility
of the wave fields, the emitted energy backpropagates
through the same environment, practically retracing the
paths of the original signal, and refocuses, again as a short
pulse, on the location
where the source emitted the original pulse. This, certainly, is
a slightly oversimplified description of the  rather complex 
physics which is involved in the real 
time-reversal experiment. In 
practice, the quality of the refocused pulse depends on
many parameters, as for example the randomness of the
background medium, the size and location of the 
receiver array, temporal fluctuations of the environment,
etc. Surprisingly, the quality of the refocused 
signal increases with increasing complexity of the 
background environment. This was observed experimentally
by M. Fink and his group at the Laboratoire Ondes et Acoustique
at Universit\'{e} Paris VII in Paris
by a series of  laboratory experiments (see
for example \cite{Draeger97,Fink97,Fink01}),
 and by W.~A.~Kuperman and co-workers at the Scripps Institution of 
Oceanography  at University of California, San Diego,
 by  a series of experiments performed
between 1996 and 2000 in a shallow ocean environment
(see for example \cite{Collins91,Edelmann02,Hodgkiss99}).

The list of possible applications of this time-reversal technique
is long. In an iterated fashion, resembling the power 
method for finding maximal eigenvalues of a square
matrix, the time-reversal method can be applied for focusing 
energy created by an ultrasound transducer array
on strongly scattering objects in the region of interest. 
This can be used for example in lithotripsy for localizing
and destroying gall-stones in an automatic way, or more
generally in the application of  medical imaging problems. 
Detection of cracks in the aeronautic industry, or of
submarines in the ocean are other examples. See for example
\cite{Prada96,Prada02}, and for related work \cite{Cheney01,Devaney}. 
The general
idea of time-reversibility, and its use in imaging and
detection, is certainly not that new. Looking into the
literature for example of seismic imaging, the application
of this basic idea can be found in a classical and very
successful imaging strategy for detecting scattering
interfaces  in the earth, the so-called 'migration' 
technique \cite{Bleistein01,Claerbout85}.
However, the systematic use and investigation 
of the time-reversal  phenomenon and its
experimental realizations 
 started more recently, and has been carried out 
during the last 10--15 years or so by different
research groups. See for example 
\cite{Draeger97,
Edelmann02,Fink97,Fink01,Hodgkiss99,Kuperman98,Song03,Yang03}
for experimental demonstrations, and 
\cite{Bal03,Bardos02,Berryman02,Blomgren02,Borcea03,
Collins91,Dowling92,Dowling93,Fouque04,Jackson91,
Klibanov03,Papanicolaou04,Prada96, Smith99, Snieder98,
Solna02, Tsogka02,Zhao04} for theoretical 
and numerical approaches.

One very young and promising application of time-reversal
is {\bf communication}. In this paper we will mainly concentrate
on that application, although the general results should
carry over also to other applications as mentioned above. 

The paper is organized as follows. In section \ref{TRO.Sec}
we give a very short introduction
 into time-reversal in the ocean, with a special
emphasis on underwater sound communication. Wireless
communication in a MIMO setup, our second 
main application in this paper,
is briefly presented in section \ref{MIMO.Sec}.
 In section
\ref{Sym.Sec} we discuss symmetric hyperbolic systems in 
the form needed here, and examples of such systems are given
in section \ref{Exa.Sec}. In section \ref{Spa.Sec}, the 
basic spaces and operators necessary for our mathematical
treatment are introduced. The inverse problem in communication,
which we are focusing on in this paper, is then defined
in section \ref{Opt.Sec}. In section \ref{Mat.Sec}, we derive
the basic iterative scheme for solving this inverse problem.
Section \ref{Adj.Sec} gives practical expressions for 
 calculating the adjoint communication
operator, which plays a key role in the iterative time-reversal
 schemes  presented in this paper.
In section \ref{ATRM.Sec} the acoustic time-reversal mirror
is defined, which will provide the link between the 
'acoustic time-reversal experiment' and the adjoint 
communication operator. The analogous results for the electromagnetic
time-reversal mirror are discussed in section \ref{ETRM.Sec}.
Section \ref{TRA.Sec} combines the results of these two sections,
and explicitly provides the link between time-reversal and
the adjoint imaging method. Sections \ref{GMC.Sec},
\ref{MNS.Sec}, and \ref{RLS.Sec} propose then
several different iterative time-reversal schemes for solving
the inverse problem of communication, using this key
relationship between time-reversal and the adjoint communication
operator. The practically important issue of partial measurements
(and generalized measurements)
is treated in section \ref{PM.Sec}. Finally, section 
\ref{SFR.Sec} summarizes the results of this paper, and
points out some interesting future research directions.

\section{Time-reversal and communication in the ocean}\label{TRO.Sec}

The ocean is a complex wave-guide for sound \cite{Collins94}.
 In addition to
scattering effects at the top and the bottom of the ocean,
also its temperature profile and the corresponding 
refractive effects contribute to this wave-guiding property and
 allow
acoustic energy to travel large distances. Typically, this
propagation has a very complicated behaviour. For example, 
multipathing occurs if source and receiver of sound waves
are far away from each other, since due to scattering and
refraction there are many 
possible  paths on which acoustic energy can travel
between them.  
Surface waves and air bubbles at the 
top of the ocean, sound propagation through
the rocks and sedimentary layers at the
bottom of the ocean, and other effects 
further contribute to the complexity of sound propagation 
in the ocean. 
When a source (e.g. the base station of a communication system
in the ocean) 
emits a short signal at a given location, the
receiver (e.g. a user of this communication system) 
some distance away from the source typically receives a long and
complicated signal due to the various influences of the ocean 
to this signal along the different connecting paths. 
If the base station wants to communicate with the
user by sending a series of short signals, this complex response
of the ocean to the signal needs to be resolved and taken into
account.

In a classical communication system, the base station
which wants to communicate with a user broadcasts a series
of short signals (e.g.~a series of 'zeros' and
'ones') into the environment. The hope is that the user 
will receive this message as a similar series of 
temporally well-resolved  short signals which 
can be easily identified  and decoded. However, this  
almost never occurs in a complex environment, 
due to the multipathing and the 
resulting delay-spread of the emitted signals. Typically,
when the base station broadcasts a series of short signals into
such an environment,
intersymbol interference  occurs at the user position 
due to the temporal overlap of the multipath contributions
of these signals. 
In order to recover the individual signals, a significant amount of 
signal processing is  necessary at the user side, and, most importantly,
the user needs to have some knowledge of
 the propagation behaviour of the signals in this environment
(i.e. he needs to know the 'channel'). Intersymbol interference 
can in principle be avoided by adding a sufficient delay between
individual signals emitted at the base station which takes
into account the delay-spread in the medium. That, however,
slows down the communication, and reduces the
capacity of the environment as a communication system. 
An additional drawback of simply  broadcasting communication signals 
from the base station is the obvious lack of interception security.  A
different user who also knows the propagation behaviour of signals 
in the environment,
 can equally well resolve the series of signals arriving 
at his location and decode them. 

Several approaches have been suggested to circumvent the above
mentioned drawbacks of communication in multiple-scattering
environments. 
Some very promising techniques are based on the time-reversibility 
property of propagating wave-fields 
\cite{Edelmann02,Emami04,Kim03,Kyritsi04,Montaldo04a,Montaldo04b,Oestges04,Smith99}.
The basic idea is as follows. The user who wants to communicate
with the base station, starts the communication process 
by sending a short pilot signal through
the environment. The base station receives this signal as a long
and complex signal due to multipathing. 
It time-reverses the received signal and sends it
back into the environment. The backpropagating waves will produce a
complicated wave-field everywhere due to the many interfering
parts of the emitted signal. However, due to  the 
 time-reversibility of the wave-fields, we
expect that the interference will be constructive at
the position of the user who sent the original pilot signal,
and mainly destructive at all other positions. Therefore,
the user will receive at his location a short signal very similar
to (ideally, a time-reversed replica of) 
the pilot signal which he sent for starting the communication.
All other users who might be in the environment at the same
time will only
receive noise speckle due to incoherently interfering contributions
of the backpropagating field. If the base station sends the
individual elements ('ones' and 'zeros') of the intended message
 in a phase-encoded form as a long overlapping string of signals,
the superposition principle will ensure that, at the user position,
this string of signals will appear as a series of short well-separated
signals, each resembling some phase-shifted (and time-reversed) 
form of the pilot signal. 

\begin{figure} 
\begin{center} 
{\hspace*{1cm}} 
\epsfysize6cm 
\epsfbox{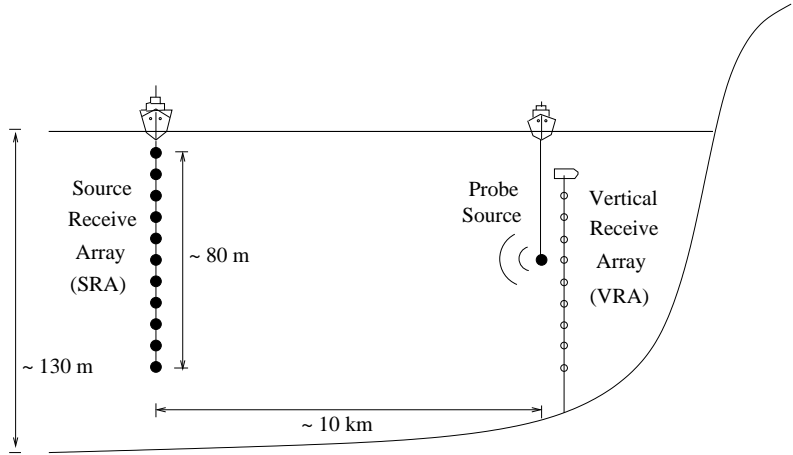} 
\vspace{1mm} 
{\small 
\caption{\label{Figure1} 
The communication problem in underwater acoustics } 
} 
\end{center} 
\end{figure}

In order to find out whether this theoretically predicted scenario actually
takes place in a real multiple-scattering environment like the 
ocean, Kuperman  \etal have performed
a series of four experiments in a shallow ocean environment between
1996 and 2000, essentially following the above described scenario.  
 The experiments have been performed at a Mediterranean
location close to the Italian coast. (A similar setup was 
also used in an  experiment performed
in 2001 off the coast of New England which 
has been reported in Yang \cite{Yang03}.) 
A schematic view of these experiments
is shown in figure \ref{Figure1}. The single 'user' is replaced  here
by a 'probe source', and the 
'source-receive array' (SRA) plays the role of the 
'base station'. An additional  vertical receive array (VRA) was
deployed at the position of the probe source in order to measure
the temporal and spatial spread of the backpropagating fields
in the neighbourhood of the probe source location. 
In this shallow environment 
(depths of about $100-200$ m, and distances between $10$ and $30$ km)
the multipathing of the waves is mostly caused by multiple reflections
at the surface and the bottom of the ocean. The results of the 
experiments have been reported in \cite{Edelmann02,Hodgkiss99,Kuperman98}. 
They show that in fact
a strong spatial and temporal focusing of the backpropagating waves
occurs at the source position. A theoretical explanation of the
temporal and spatial refocusing of time-reversed waves in 
random environments has been  given in Blomgren \etal \cite{Blomgren02}.

\section{The MIMO setup in wireless communication}\label{MIMO.Sec}
\begin{figure} 
\begin{center} 
{\hspace*{1cm}} 
\epsfysize6cm 
\epsfbox{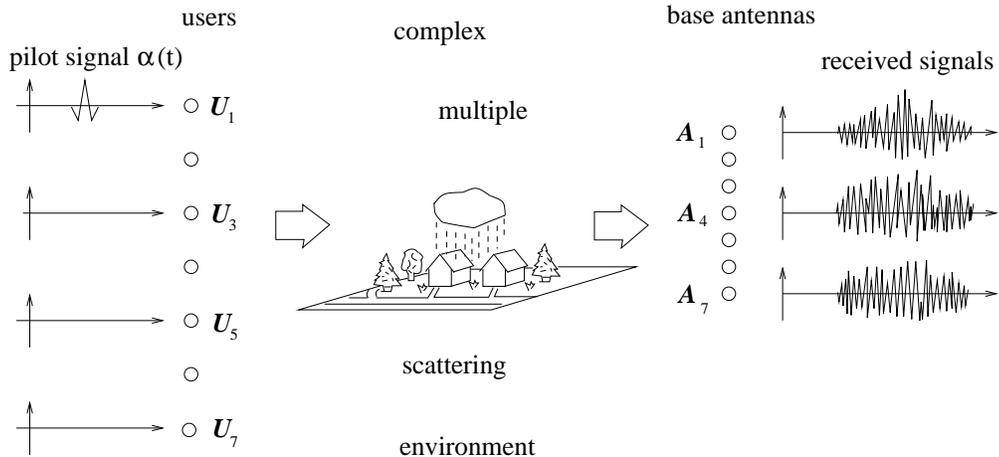} 
\vspace{1mm} 
{\small 
\caption{\label{Figure2} 
Schematic view of a typical MIMO setup in wireless communication. } 
} 
\end{center} 
\end{figure}

The underwater acoustics scenario described above 
directly carries over to situations which might be
more familiar to most of us, namely to the more and
more popular wireless 
communication networks using mainly 
electromagnetic waves in the microwave regime. 
Starting from the everyday use of cell-phones, 
ranging to small wireless-operating local area networks (LAN)
for computer systems or for private enterprise
communication systems, exactly the same problems arise
as in underwater communication. The typically employed
microwaves of a wavelength at about 10--30 cm are heavily scattered
by  environmental objects like cars, fences, trees, doors, 
furniture etc. This causes a very complicated 
multipath structure of the signals received by users of 
such a communication system. Since bandwidths are limited
and increasingly expensive, a need for more and more efficient 
communication systems is imminent. Recently, the idea of 
a so-called multiple-input multiple-output (MIMO) communication
system has been introduced with the potential to increase
the capacity and efficiency of wireless communication systems
\cite{Foschini96}. 
The idea is to replace a single antenna at the base station
which is responsible for multiple users, 
or even a system where
one user communicates only with his own dedicated base antenna, 
 by a more general system  where an array of multiple antennas
at the base station is interacting simultaneously and in a 
complex way 
with multiple users. A schematic description of such a MIMO system
(with seven base antennas and seven users) is given in figure \ref{Figure2}.
See for example \cite{Gesbert03,Simon01} for recent overviews on MIMO
technology.

Time-reversal techniques are likely to play also here a key
role in improving communication procedures and for optimizing
the use of the limited resources (especially bandwidth)
which are available for this technology
\cite{Emami04,Kim03,Kyritsi04,Oestges04}.
One big advantage of time-reversal techniques is that
they are automatically adapted to the complex 
environment and that they can be very fast since
they do not require heavy signal processing at the receiver
or user side. 
In \cite{Montaldo04a,Montaldo04b},  an  {\bf iterated 
time-reversal scheme} for the optimal refocusing of signals
in such a MIMO communication system was proposed, which we will describe
in more details in section \ref{GMC.Sec}.

In the present paper we will establish a direct link between 
the time-reversal technique and solution strategies
for inverse problems.
 As an application of this relationship,
we will derive iterative time-reversal schemes
for the optimization of 
wireless or underwater acoustic MIMO communication
systems.
The derivation is performed completely in
time-domain, for very general first order symmetric hyperbolic
systems describing wave propagation phenomena in a complex
environment. 
One of the schemes which we derive, in a certain
sense the 'basic one', will turn out to be 
practically equivalent to the scheme introduced 
in  \cite{Montaldo04a,Montaldo04b}, although the
derivation uses different tools.  Therefore, it provides a 
new interpretation of that scheme. The other schemes
which we introduce are new in this application, and can be considered
as either generalizations of the basic scheme,
or as independent alternatives to that scheme. 
Each of them addresses slightly different
objectives and has its own very specific 
characteristics.

\section{Symmetric hyperbolic systems}\label{Sym.Sec} 
We  treat wave propagation in communication systems
in the general framework of 
 symmetric hyperbolic systems of the form 
\begin{equation} 
    \label{Sym.1}
\Gamma(\bfx) \frac{\partial \bfu}{\partial t}\,+\,\sum_{i=1}^{3}D^{i}  
\frac{\partial \bfu}{\partial x_i} \,+\, \Phi (\bfx)  \bfu
\,=\, \bfq 
\end{equation} 
\begin{equation} 
   \label{Sym.2}
\bfu(\bfx,0) \,=\, 0. 
\end{equation} 
Here, $\bfu(\bfx,t)$ and $\bfq(\bfx,t)$ 
are real-valued time-dependent $N$-vectors, 
$\bfx\in\R^3$, and 
$t\in [0,T]$. 
$\Gamma(\bfx)$ is a real, symmetric, uniformly positive 
definite $N\times N$-matrix, i.e., $\Gamma(\bfx)\geq \epsilon$ for some  
$\epsilon > 0$. 
Moreover, $\Phi(\bfx)$ is a symmetric positive semi-definite $N\times N$
matrix, i.e., $\Phi(\bfx)\geq 0$. It models possible energy loss through 
dissipation in the medium.
The $D^i$ are real,  symmetric and  
independent of $(\bfx,t)$.  
We will also use the short notation
\begin{equation} 
   \label{Sym.3} 
\Lambda\,:=\, \sum_{i=1}^{3}D^i  
\frac{\partial}{\partial x^i}.
\end{equation} 
In addition to the above mentioned assumptions on the coefficients
$\Gamma(\bfx)$ and $\Phi(\bfx)$, we will assume throughout this paper 
that all quantities $\Gamma(\bfx)$, $\Phi(\bfx)$, $\bfu(\bfx,t)$ and
 $\bfq(\bfx,t)$ are
'sufficiently regular' in order to safely apply for example integration
by parts and Green's formulas. For details see for example
\cite{Dautray,Marsden}.
 Throughout this paper we 
will assume that no energy
reaches the boundaries
$\partial \Omega$ during the time $[0,T]$, such that we will
always have 
\begin{equation}\label{Sym.4}
\bfu(\bfx,t)=0\quad \mbox{on} \quad\partial\Omega\times [0,T].
\end{equation}
The {\sl energy density} $\mcalE(\bfx,t)$ is defined by 
\begin{eqnarray*} 
\mcalE(\bfx,t)&=& \frac{1}{2}\left\langle \Gamma(\bfx)\bfu(\bfx,t), 
\bfu(\bfx,t)\right\rangle_{N} \\ 
  &=&
\frac{1}{2} \sum_{m,n=1}^{N}\Gamma_{mn}(\bfx)\bfu_m(\bfx,t)\bfu_n(\bfx,t).
\end{eqnarray*} 
The {\sl total energy} $\hcalE(t)$ in $\Omega$ 
at a given time $t$ is therefore
\begin{displaymath} 
\hcalE(t)\,=\,\frac{1}{2}\int_{\Omega} \left\langle \Gamma(\bfx)\bfu(\bfx,t), 
\bfu(\bfx,t)\right\rangle_{N}\,d\bfx.
\end{displaymath}
The {\sl flux} $\mcalF(\bfx,t)$ is given by 
\begin{displaymath} 
\mcalF_i(\bfx,t)\,=\,\frac{1}{2}\left\langle D^i\bfu(\bfx,t),
\bfu(\bfx,t)\right\rangle_{N}\quad,\quad i=1,2,3.
\end{displaymath}

\section{Examples for symmetric hyperbolic systems}\label{Exa.Sec} 
In the following, we want to give some examples for symmetric 
hyperbolic systems as defined above. Of special interest for
communication are the system of the linearized acoustic equations
and the system of Maxwell's equations. We will discuss these two
examples in detail in this paper.
Another important example for symmetric
hyperbolic systems is the system of elastic waves equations, which we
will however leave out in our discussion for the sake of shortness. We only
mention here that all the results derived here apply without restrictions
also to linear elastic waves. Elastic wave propagation becomes important
for example in ocean acoustic communication models which incorporate
wave propagation through the sedimentary and rock layers at the
bottom of the ocean. 

\subsection{Linearized acoustic equations} \label{Exa.Aco}
As a model for {\bf underwater sound propagation}, we consider the following 
linearized form of the acoustic equations in an isotropic medium
\begin{eqnarray}
   \label{Exa.1.1}
\rho(\bfx) \frac{\partial \bfv}{\partial t}\,+\, \ngrad \, p 
 &=& \bfq_{v} \\
   \label{Exa.1.2} 
\kappa(\bfx) \frac{\partial p}{\partial t}\,+\, \ndiv \, \bfv
 &=& \bfq_{p} 
\end{eqnarray} 
\begin{equation}
   \label{Exa.1.3} 
p(\bfx,0)=0,\qquad \bfv(\bfx,0)=0.  
\end{equation} 
Here, $\bfv$ is the velocity, $p$ the pressure, $\rho$ the density, and 
$\kappa$ the compressibility.  
We have $N=4$, $\bfu=(\bfv,p)^{T}$ and $\bfq=(\bfq_v,\bfq_p)^{T}$
(where '$T$' as a superscript always means 'transpose'). Moreover, we have
\begin{displaymath} 
\Gamma(\bfx)\,=\,\ndiag(\rho(x),\rho(x),\rho(x),\kappa(x))
\quad \mbox{and} \quad \Phi(\bfx)\,=\,0. 
\end{displaymath}  
With the notation $\varphi=(\partial_1,\partial_2,\partial_3)^T$, 
we can write $\Lambda$ as 
\begin{displaymath}
\Lambda\,=\,\left( \begin{array}{cc}
0 & \varphi \\ \varphi^T & 0
\end{array}\right).
\end{displaymath} 
The operators $D^i$, $i=1,2,3$, can be recovered here from $\Lambda$ by 
putting 
\begin{displaymath}
D_{m,n}^{i}\,=\,\left\{
\begin{array}{cl}
1 &\mbox{where}\quad\Lambda_{m,n} = \partial_i,\\
0 & \mbox{elsewhere}
\end{array}\right.
\end{displaymath}
The energy density $\mcalE(\bfx,t)$ is given by 
\begin{displaymath}
\mcalE(\bfx,t)\,=\, \frac{1}{2}\big( \rho(\bfx)|\bfv(\bfx,t)|^2\,+\,
 \kappa(\bfx) p^2(\bfx,t)\big),
\end{displaymath} 
and the energy flux $\mcalF(\bfx,t)$ is
\begin{displaymath}
\mcalF(\bfx,t)\,=\,p(\bfx,t)\bfv(\bfx,t).
\end{displaymath} 
We mention that the dissipative case $\Phi(\bfx)\neq 0$ can be treated as well
in our framework, and yields analogous results to those presented
here. 
 
\subsection{Maxwell's equations}\label{Exa.Ele} 
As a second example, we will consider Maxwell's equations for an
anisotropic medium with some energy loss due to the
inherent conductivity.  This can for example model 
{\bf wireless communication} in a complex  environment.
\begin{eqnarray} 
   \label{Exa.2.1} 
\epsilon(\bfx) \frac{\partial \bfE}{\partial t}\,-\, \ncurl \,\bfH  
\,+\,\sigma(\bfx) \bfE &=& \bfq_{E} \\ 
   \label{Exa.2.2} 
\mu(\bfx) \frac{\partial \bfH}{\partial t}\,+\, \ncurl \,\bfE  
\qquad \qquad &=& \bfq_{H}
\end{eqnarray} 
\begin{equation}
   \label{Exa.2.3} 
\bfE(\bfx,0)=0,\qquad \bfH(\bfx,0)=0.  
\end{equation} 
We have $N=6$, and $\bfu=(\bfE,\bfH)^T$, $\bfq=(\bfq_E,\bfq_H)^T$. 
Moreover,  
\begin{eqnarray*}
\Gamma(\bfx)&=&\ndiag(\epsilon(\bfx),\mu(\bfx))\\
\Phi(\bfx)&=&\ndiag(\sigma(\bfx),0).
\end{eqnarray*} 
Here, $\epsilon$ and $\mu$ are symmetric positive definite 
$3\times 3$ matrices modelling the anisotropic permittivity and
permeability distribution in the medium, and $\sigma$ is
a symmetric positive semi-definite $3\times 3$ matrix which 
models the anisotropic conductivity distribution.
In wireless communication, this form can model for example 
dissipation by conductive trees, 
conductive wires, rainfall, pipes, etc.
The operator $\Lambda$ can be written in block form as 
\begin{displaymath}
\Lambda\,=\,  \left( \begin{array}{cc}
0 & -\Xi \\
\Xi & 0
\end{array}\right) ,
\end{displaymath}
with
\begin{displaymath}
\Xi\,=\, \left( \begin{array}{ccc}
0 &- \partial_3 & \partial_2 \\
\partial_3 & 0 & -\partial_1 \\
-\partial_2 & \partial_1 & 0 
\end{array}\right) . 
\end{displaymath}
The operators $D^i$, $i=1,2,3$, can be recovered here from $\Lambda$ by 
putting 
\begin{displaymath}
D_{m,n}^{i}\,=\,\left\{
\begin{array}{cl}
1 &\mbox{where}\quad\Lambda_{m,n} = \partial_i,\\
-1 & \mbox{where}\quad\Lambda_{m,n}= -
\partial_i,\\
0 & \mbox{elsewhere}
\end{array}\right.
\end{displaymath}
The energy density  $\mcalE(\bfx,t)$ is given by 
\begin{displaymath}
\mcalE(\bfx,t)\,=\,\frac{1}{2}\big( \epsilon(\bfx)|\bfE(\bfx,t)|^2\,+\,
\mu(\bfx)|\bfH(\bfx,t)|^2 \big).
\end{displaymath} 
The energy flux $\mcalF(\bfx,t)$ is 
described by the {\sl Poynting vector}
\begin{displaymath}
\mcalF(\bfx,t)\,=\,\bfE(\bfx,t)\times\bfH(\bfx,t).
\end{displaymath}

\subsection{Elastic waves equations}

As already mentioned, also elastic waves can be treated in the
general framework of symmetric hyperbolic systems. For more details
 we refer to \cite{Marsden,Ryzhik96}.

\section{The basic spaces and operators}\label{Spa.Sec}
 
For our mathematical treatment of time-reversal
 we introduce the following 
spaces and inner products.
We assume that we have $\anzus$ users 
 $\mcalU_j$, $j=1,\ldots,\anzus$ in our system, and 
in addition a base station which consists of 
$\anzant$  antennas
$\mcalA_k$, $k=1,\ldots,\anzant$. Each user and each antenna
at the base station can receive, process and emit signals which 
we denote by $\bfs_j(t)$ for a given user 
$\mcalU_j$, $j=1,\ldots,\anzus$, and by 
$\bfr_k(t)$ for a given base antenna
$\mcalA_k$, $k=1,\ldots,\anzant$. 
Each of these signals consists of a time-dependent $N$-vector
indicating measured or processed signals of time-length $T$.
In our analysis we will often have to consider functions defined
on the time interval $[T,2T]$ instead of $[0,T]$. 
For simplicity, we will use the same notation for the function
spaces  defined on  $[T,2T]$ as we use for those defined
on $[0,T]$. It will always be obvious which space we refer
to in a given situation. 

Lumping together all signals  at
the users on the one hand, and all signals at the base station
on the other hand, yields the two fundamental quantities 
\begin{eqnarray*}
\bfs \,=\, (\bfs_1,\ldots,\bfs_{\anzus}) &\in&  \Zhat \\
\bfr \,=\, (\bfr_1,\ldots,\bfr_{\anzant}) &\in&  Z 
\end{eqnarray*}
with 
\begin{displaymath}
\Zhat\,=\,\big(L_2([0,T])^N\big)^{\anzus},
\qquad Z\,=\,\big(L_2([0,T])^N\big)^{\anzant}.
\end{displaymath}
The two signal spaces $Z$ and $\Zhat$ introduced above 
are equipped with the inner products
\begin{eqnarray*}
\left\langle\bfs^{(1)},\bfs^{(2)}\right\rangle_{\Zhat}&=&
\sum_{j=1}^{\anzus}\int_{[0,T]}\Big\langle \bfs^{(1)}_j(t),
\bfs^{(2)}_j(t)\Big\rangle_N\,dt \\
\left\langle\bfr^{(1)},\bfr^{(2)}\right\rangle_{Z}&=&
\sum_{k=1}^{\anzant}\int_{[0,T]}\Big\langle \bfr^{(1)}_k(t),
\bfr^{(2)}_k(t)\Big\rangle_N\,dt.
\end{eqnarray*}
The corresponding norms are 
\begin{displaymath}
   \label{4.2}
\|\bfs\|_{\Zhat}^2\,=\,\left\langle\bfs,\bfs\right\rangle_{\Zhat},\qquad
\|\bfr\|_{Z}^2\,=\,\left\langle\bfr,\bfr\right\rangle_{Z}
\end{displaymath}

Each user $\mcalU_j$ and each antenna $\mcalA_k$
at the base station can send a given
signal $\bfs_j(t)$ or $\bfr_k(t)$, respectively. This gives rise to 
a source distribution $\bfqh_j(\bfx,t)$, $j=1,\ldots,\anzus$,
or $\bfq_k(\bfx,t)$, $k=1,\ldots,\anzant$, respectively. Here
and in the following we will use in our
notation the following convention. If one symbol
appears in both forms, with and without a 'hat' (\^{})  on top of this
symbol, 
then all quantities {\sl with} the 'hat' symbol are related to the users, 
and those {\sl without} the 'hat' symbol to the antennas at the base station. 

Each of the sources created by a user or by a base antenna will
appear on the right hand side of (\ref{Sym.1}) as a mathematical
source function and gives rise to a  corresponding wave field
which satisfies (\ref{Sym.1}), (\ref{Sym.2}). When solving the system
(\ref{Sym.1}), (\ref{Sym.2}), typically certain Sobolev spaces 
need to be employed for the appropriate description of the underlying 
function spaces (see for example \cite{Dautray}). 
For our purposes, however, it will be sufficient
to assume that both, source functions and wave fields,
are members of the following canonical function space $U$ which 
is defined as 
\begin{displaymath}
U=\{\bfu\in L^2(\Omega\times [0,T])^{N}\,,\,\bfu=0\;
\mbox{on}\;\partial \Omega\times [0,T]\,,\,\;\|\bfu\|_{U}<\infty\},
\end{displaymath}
and which we have equipped with the usual 
 energy inner product 
\begin{displaymath}
\left\langle\bfu(t),\bfv(t)\right\rangle_{U}\,=\,
\int_{[0,T]}\int_{\Omega}\left\langle\Gamma(\bfx)\bfu(\bfx,t),
\bfv(\bfx,t)\right\rangle_N\,d\bfx dt,
\end{displaymath}
and the  corresponding energy norm  
\begin{displaymath}
\|\bfu\|_{U}^2\,=\,\left\langle\bfu\,,\,\bfu\right\rangle_{U}.
\end{displaymath}
Also here, in order to simplify the notation, we will use the
same space  when considering functions in the shifted time
interval $[T,2T]$ instead of $[0,T]$.

Typically, when a user or an antenna at the base station 
transforms a signal into a source distribution, it is done
according to a very specific antenna characteristic which
takes into account the spatial extension of the user or
the antenna. 
We will model this characteristic at the user by the functions
$\gammahat_j(\bfx)$, $j=1,\ldots,\anzus$,
and for base antennas by the functions
$\gamma_k(\bfx)$,  $k=1,\ldots,\anzant$. With these functions,
we can introduce the linear 'source operators' $\Qhat$ and $Q$ mapping
signals $\bfs$ at the set of users and $\bfr$ at the set of 
base antennas into the corresponding source
distributions $\bfqh(\bfx,t)$ and $\bfq(\bfx,t)$, respectively. 
They are given as 
\begin{eqnarray*}
  \Qhat\,:\,\Zhat \rightarrow U, &\quad& \bfqh(\bfx,t)\,=\,\Qhat \bfs
\,=\, \sum_{j=1}^{\anzus}\gammahat_j(\bfx)\bfs_j(t), \\
 Q\,:\,Z \rightarrow U, &\quad& \bfq(\bfx,t)\,=\,Q \bfr
\,=\, \sum_{k=1}^{\anzant}\gamma_k(\bfx)\bfr_k(t).
\end{eqnarray*}
We will assume that the functions $\gammahat_j(\bfx)$ are supported 
on a small neighbourhood $\Vhat_j$ of the user location $\hat{\bfd}_j$, and 
that the functions $\gamma_k(\bfx)$ are supported on a small
neighbourhood $V_k$ of the antenna location $\bfd_k$. Moreover,
all these neighbourhoods are strictly disjoint to each other. For example,
the functions $\gammahat_j(\bfx)$ could be assumed to be $L_2$-approximations
of the Dirac delta measure $\delta(\bfx-\hat{\bfd}_j)$ concentrated at
the user locations $\hat{\bfd}_j$, and the functions $\gamma_k(\bfx)$ could 
 be assumed to be $L_2$-approximations
of the Dirac delta measure $\delta(\bfx-\bfd_k)$ concentrated at
the antenna  locations $\bfd_k$.

Both, users and base antennas can also record incoming
fields $\bfu\in U$ and transform the recorded information into signals. 
Also here, this is usually done according to very specific antenna 
characteristics of each user and each base antenna. For simplicity
(and without loss of generality),
we will assume that the antenna characteristic of a user or base antenna
for receiving signals is the same as for transmitting signals, namely
$\gammahat_j(\bfx)$ for the user and $\gamma_k(\bfx)$ for a base antenna.
(The case of more general source and measurement operators is
 discussed in  section \ref{PM.Sec}.) 
With this, we can define the linear 'measurement operators'
 $\Mhat\,:\,U \rightarrow \Zhat$ and
$M\,:\,U \rightarrow Z$, respectively, which transform incoming 
fields into measured signals, by 
\begin{eqnarray*}
  \bfs_j(t)&=&(\Mhat \bfu)_j
\,=\,\int_{\Omega} \gammahat_j(\bfx)\bfu(\bfx,t)d\bfx, \quad (j=1,\ldots, \anzus) \\
 \bfr_k(t)&=&(M \bfu)_k
\,=\, \int_{\Omega}\gamma_k(\bfx)\bfu(\bfx,t)d\bfx,  \quad (k=1,\ldots, \anzant) 
\end{eqnarray*}
Finally, we define the linear operator $F$ mapping sources $\bfq$ to states 
$\bfu$ by
\begin{displaymath}
F:\,U\rightarrow U,\qquad F\bfq\,=\,\bfu,
\end{displaymath} 
where $\bfu$ solves the problem (\ref{Sym.1}), (\ref{Sym.2}). 
As already mentioned, we assume that the domain $\Omega$ is chosen 
sufficiently large
and that the boundary $\partial\Omega$ is sufficiently
far away from the users and  base antennas, such that
there is no energy reaching the boundary
in the time interval $[0,T]$ (or $[T,2T]$) due to the finite speed of signal
propagation. Therefore, the operator $F$ is well-defined.

Formally, we can now introduce the two linear {\bf communication operators}
$A$ and $B$ which are at the main focus of this paper. 
They are defined as 
\begin{eqnarray*}
A:\,Z\rightarrow \Zhat, &\quad& A\bfr\,=\,\Mhat FQ\bfr, \\
B:\,\Zhat\rightarrow Z, &\quad& B\bfs\,=\,M F \Qhat \bfs.
\end{eqnarray*}
The operator $A$ models the following situation. The base
station emits the signal $\bfr(t)$ which propagates through
the complex environment. The users measure the arriving 
wave fields and transform them into measurement signals.
The measured signals at the set of 
all users is  $\bfs(t)=A\bfr(t)$.  The operator
$B$ describes exactly the reversed situation. All users
emit together the set of signals $\bfs(t)$, which propagate
through the given complex environment and are received
by the base station. The corresponding set of measured signals at
all antennas of 
the base station is just $\bfr(t)\,=\,B\bfs(t)$. No time-reversal
is involved so far.

\section{An inverse problem arising in communication}\label{Opt.Sec} 
 
In the following, we outline a typical problem arising
in communication, which gives rise to a mathematically
well-defined inverse problem. 

A specified user of the system, say $\mcalU_1$, defines a (typically
but not necessarily short) 
pilot signal $\alpha(t)$ which 
he wants to use as a template for receiving the information from the base
station.  The base station wants to emit a
signal $\tilde{\bfr}(t)$ which, after having travelled
through the complex environment and arriving
at the user $\mcalU_1$,  matches this pilot signal
as closely as possible. Neither the base station nor
any other user except of $\mcalU_1$  
are required (or expected) to know the correct
form of the pilot signal $\alpha(t)$ for this problem.
 As an additional constraint, 
the base station wants that at the other users $\mcalU_j$, $j>1$,
as little energy as possible arrives when communicating
with the specified user  $\mcalU_1$. This is also in the
interest of the other users, who want to use a different
'channel' for communicating at the same time with the
base antenna, and want to minimize interference with the
communication initiated by user $\mcalU_1$. The complex environment
itself in which  the communication 
takes place (i.e. the 'channel') is assumed to be unknown to all
users and to the base station.

In order to arrive at a mathematical description
of this problem, we define the 'ideal signal' $\tilde{\bfs}(t)$
received by all users as
\begin{equation}\label{Mat.1}
\tilde{\bfs}(t)\,=\,(\alpha(t),0,\ldots,0)^T.
\end{equation}
Each user only knows his own component of this signal,
and the base antenna does not need to know any component
of this ideal signal at all. 
\begin{definition}
{\rm {\bf The inverse problem of communication}:\quad
In the terminology of inverse problems, the above described scenario
defines an inverse source problem,
which we call for the purpose of this paper the 
 'inverse problem of communication'. 
The goal is to 
find a 'source distribution' $\rtilde(t)$ at the base 
station which satisfies the 'data' $\stilde(t)$ at the users:
\begin{equation}\label{Opt.2}
A\rtilde\,=\,\stilde.
\end{equation}
The 'state equation' relating sources to data is given
by the symmetric hyperbolic system (\ref{Sym.1}), (\ref{Sym.2}).
}
\end{definition}

\begin{remark}
{\rm
Notice that the basic operator $A$ in (\ref{Opt.2}) is unknown to
us since we do not know the complicated medium
in which the waves propagate.
If the operator $A$ (together with $\stilde$)
would be known at the base station
by some means, the inverse
source problem formulated above could be solved using classical
inverse problems techniques, which would be computationally
expensive but in principle doable. In our situation, we are 
able to do physical experiments, which amounts to 'applying' the
communication operator $A$ to a given signal. Determining
the operator $A$ explicitly by applying it to a set of basis functions
of $Z$ would be possible, but again too expensive. 
We will show in the following that, nevertheless, many of the
classical solution schemes known from inverse problems theory
can be applied in our situation even without knowing this
operator. The basic tool which we will use is a series
of time-reversal experiments, applied to carefully designed
signals at the users and the base station.
} 
\end{remark}

\begin{remark}
{\rm 
A practical template for an {\bf iterative scheme for finding an optimal
signal} at the base station can be envisioned as 
follows. User  $\mcalU_1$
starts the communication 
process  by emitting an initial signal $\bfs_1^{(0)}(t)$
into the complex environment. This signal, after having
propagated through  the complex environment, 
finally arrives at the base station and is received there usually
as a relatively long and complicated signal due to the 
multiple scattering events it experienced on its way.
  When the base station
receives such a signal, it processes it and sends a 
new signal $\bfr^{(1)}(t)$ back through 
the environment
which is received by all users. After receiving this 
signal, all users become active. The user $\mcalU_1$ compares
the received signal with the pilot signal. If the match
is not good enough, it processes
the received signal in order to optimize the match with the
pilot signal when receiving the next iterate from the base station.
All other users identify the received signal as unwanted noise, and 
process it with the goal to receive in the next iterate from the 
base station a signal with lower amplitude, such that it does
not interfere with their own communications. 
All users send now their processed signals, which together 
define $\bfs^{(1)}(t)$, 
back to the base station. The base station 
receives them all at the same time, 
again usually as a long and complicated signal, processes this
signal and
sends a new signal $\bfr^{(2)}(t)$ back into the environment
which ideally will match the desired signals at all users
better than the previously emitted signal $\bfr^{(1)}(t)$.
This iteration stops when all users are satisfied, i.e.,
when user $\mcalU_1$ receives a signal which is sufficiently
close to the pilot signal $\alpha(t)$, and the energy or
amplitude of the signals arriving at the other users has decreased
enough in order not to disturb their own communications. 
After this learning process of the channel has been completed,
 the user $\mcalU_1$ can now start communicating safely
with the base antenna using the chosen pilot signal 
$\alpha(t)$ for decoding the received signals. 
}
\end{remark}

Similar schemes have been suggested in 
\cite{Edelmann02,Emami04,Kim03,Kyritsi04,Oestges04,Smith99}
in a single-step fashion,  and
in \cite{Montaldo04a,Montaldo04b} performing multiple steps
of the iteration.
The main questions to be answered are certainly which signals
each user and each base antenna needs to emit in each step,
how these signals need to be processed, at which stage this
iteration should be terminated, and which 
optimal solution this scheme is expected to converge to. 
{\em One goal of this paper is to provide a theoretical framework
for answering these questions by combining basic
concepts of inverse problems theory with experimental 
time-reversal techniques.}

\section{The basic approach for solving the 
inverse problem of communication}\label{Mat.Sec}

A standard approach for solving  problem (\ref{Opt.2}) in the
situation of noisy data is to look for the least-squares solution 
\begin{equation} \label{Mat.2}
\bfr_{LS}\,=\, \mbox{Min}_{\bfr}\, 
\| A \bfr - \tilde{\bfs}\|_{\Zhat}^2
\end{equation} 
In order to practically find a solution of (\ref{Mat.2}), we introduce
 the cost functional 
\begin{equation} \label{Mat.3}
\mcalJ(\bfr)\,=\,\frac{1}{2}\Big\langle A\bfr - \tilde{\bfs} 
\,,\, A\bfr - \tilde{\bfs}\Big\rangle_{\Zhat}
\end{equation}
In the {\bf  basic approach} we propose to use the 
{\bf gradient method} for finding the minimum of (\ref{Mat.3}). 
In this method, in each iteration a correction $\delta \bfr$ is sought
for a guess $\bfr$ 
which points into the negative gradient direction
$-A^{*}(A\bfr-\tilde{\bfs})$ of the cost functional (\ref{Mat.3}).
In other words, starting with the initial guess 
$\bfr^{(0)}=0$, the iteration of the gradient method goes
as follows:
\begin{equation}  \label{Mat.4}
\begin{array}{rcl}
\bfr^{(0)}&=&0 \\ 
\bfr^{(n+1)}&=& \bfr^{(n)}\,-\,\beta^{(n)}A^{*}(A\bfr^{(n)}-\tilde{\bfs}),
\end{array}
\end{equation}
where $\beta^{(n)}$ is the step-size at iteration number $n$.
Notice that the signals $\bfr^{(n)}$ are measured at the 
base station, whereas the difference $A\bfr^{(n)}-\tilde{\bfs}$
is determined at the users. In particular, $\tilde{\bfs}$ is the
pilot signal only known by the user who defined it, combined with
zero signals at the remaining users. $A\bfr^{(n)}$ is the signal
received by the users at the $n$-th iteration step.

\section{The adjoint operator $A^{*}$}\label{Adj.Sec}

We see from (\ref{Mat.4}) that we will have to apply the adjoint 
operator $A^{*}$ repeatedly when implementing the
gradient method. In this section we provide practically
useful expressions for applying this operator to 
a given element of $\Zhat$.

\begin{lemma} \label{Lemma.Adj.1}
We have 
\begin{equation}
   \label{Adj.1} 
A^{*}\,=\,Q^{*}F^{*}\Mhat^{*}.
\end{equation}
\end{lemma}
\noindent {\sl Proof:} \quad  This follows from $A=\Mhat F Q$.\hfill $\Box$

\begin{theorem} \label{Theorem.Adj.1}
We have 
\begin{equation}
  \label{Adj.2}
\begin{array}{rcl} 
M^{*}\,=\,\Gamma^{-1}Q, &\quad& \Mhat^{*}\,=\,\Gamma^{-1}\Qhat,\\
Q^{*}\,=\,M \Gamma , &\quad& \Qhat^{*}\,=\,\Mhat \Gamma.
\end{array}
\end{equation}
\end{theorem}

\noindent {\sl Proof:} \quad  The proof is given in appendix A.\hfill $\Box$

Next, we want to find an expression for $F^{*}\bfv$, $\bfv\in U$, where 
$F^{*}$ is the adjoint of the operator $F$. 
\begin{theorem} \label{Theorem.Adj.2}
Let $\bfz$ be the solution of the adjoint symmetric hyperbolic system
\begin{equation}  
    \label{Adj.3}
-\Gamma(\bfx) \frac{\partial \bfz}{\partial t}\,-\,\sum_{i=1}^{3}D^{i}  
\frac{\partial \bfz}{\partial x_i}
\,+\, \Phi(\bfx)\bfz
 \,=\, \Gamma(\bfx)  \bfv(\bfx,t),
\end{equation} 
\begin{equation}  
   \label{Adj.4}
\bfz(\bfx,T) \,=\, 0 ,
\end{equation} 
\begin{equation}  
   \label{Adj.5}
\bfz(\bfx,t)=0\quad\mbox{{\rm on}}\quad\partial\Omega\times [0,T].
\end{equation} 
Then
\begin{equation}  
   \label{Adj.6}
F^{*}\bfv\,=\,\Gamma^{-1}(\bfx)\bfz(\bfx,t).
\end{equation} 
\end{theorem}
\noindent {\bf Proof:} \quad The proof is given in appendix B.\hfill $\Box$

\begin{remark}
{\rm
This procedural characterization of the adjoint operator is often used  
in solution strategies of large scale inverse 
problems, where it naturally leads to so-called 'backpropagation
strategies'. See for example 
\cite{Dorn98,Dorn99,Haber00,Hursky04,Natterer95,Vogeler03}
and the references given there.
}
\end{remark}

\begin{remark}
{\rm
Notice that in the adjoint system (\ref{Adj.3})--(\ref{Adj.5}) 
 'final value conditions' are given at $t=T$ in contrast to 
(\ref{Sym.1})--(\ref{Sym.4}) where 'initial value conditions'
are prescribed at $t=0$. This corresponds to the fact that time is 
running backward in (\ref{Adj.3})--(\ref{Adj.5}) and forward
in (\ref{Sym.1})--(\ref{Sym.4}). 
}
\end{remark}

\section{The acoustic time-reversal mirror}\label{ATRM.Sec} 

In the following we want to define an operator $\mcalS_a$ 
 such that $F^{*}=\Gamma^{-1}\mcalS_a F \mcalS_a\Gamma$ holds.
 We will call this
operator $\mcalS_a$ the '{\bf acoustic time-reversal operator}'.
We will also define the '{\bf acoustic time-reversal mirrors}'
 $\mcalT_a$ and  $\mcalThat_a$,
which act on the signals instead of the sources or fields.

We consider the acoustic system 
\begin{eqnarray} 
   \label{Atrm.1}
\rho\frac{\partial \bfv_f}{\partial t}\,+\, \ngrad   
p_f(\bfx,t) &=& \bfq_{\bfv},  \\
   \label{Atrm.2} 
\kappa\frac{\partial p_f}{\partial t}\,+\, \ndiv   
\bfv_f(\bfx,t) &=& \bfq_p, 
\end{eqnarray} 
\begin{equation}
   \label{Atrm.3} 
 \bfv_f(\bfx,0)\,=\, 0,\qquad  
p_f(\bfx,0)\,=\, 0\quad \mbox{in}\quad \Omega
\end{equation} 
with $t\in [0,T]$ and zero boundary conditions 
 at $\partial\Omega\times [0,T]$.
We want to calculate the action of the {\sl adjoint operator}
$F^{*}$
on a vector $(\phi,\psi)^T\in U$.

\begin{theorem} 
   \label{Theorem.Atrm.1}
Let $(\phi,\psi)^T\in U$ and let 
$(\bfv_a,p_a)^T$ be the solution of the 
adjoint system 
\begin{eqnarray}
   \label{Atrm.4} 
-\rho\frac{\partial \bfv_a}{\partial t}\,-\, \ngrad  
p_a(\bfx,t)  &=& \rho(\bfx) \phi(\bfx,t)  \\
   \label{Atrm.5} 
-\kappa\frac{\partial p_a}{\partial t}\,-\, \ndiv  
\bfv_a(\bfx,t) &=& \kappa(\bfx) \psi(\bfx,t) 
\end{eqnarray} 
\begin{equation}
   \label{Atrm.6} 
\bfv_{a}(\bfx,T)\,=\, 0,\qquad p_{a}(\bfx,T)\,=\,0\quad \mbox{in}
\quad \Omega,
\end{equation}
with  $t\in [0,T]$ and zero boundary conditions at
 $\partial\Omega\times [0,T]$.
Then we have 
\begin{equation}
   \label{Atrm.7} 
F^{*}\left( {\phi \atop \psi}\right)\,=\,\left( {\rho^{-1} 
\bfv_a(\bfx,t) \atop \kappa^{-1} p_a(\bfx,t)}\right).
\end{equation} 
\end{theorem}
\noindent {\sl Proof:}\quad This theorem is just an application of 
theorem \ref{Theorem.Adj.2} to the acoustic symmetric hyperbolic
system. For the convenience of the reader, we will give a direct
proof as well in appendix C.

\begin{definition}
We define the {\bf acoustic time-reversal operator} $\mcalS_a$ by 
putting for all $\bfq=(\bfq_{\bfv},\bfq_p)^{T}\in U$
(and similarly for all $\bfu=(\bfv,p)^{T}\in U$)
\begin{equation}
   \label{Atrm.8}
\big(\mcalS_a \bfq\big)(\bfx,t)\,=\, 
\left( {-\bfq_{\bfv}(\bfx,2T-t) \atop \bfq_p(\bfx,2T-t)}\right)
\end{equation}
\end{definition}
\begin{definition}
We define the {\bf acoustic time-reversal mirrors} $\mcalT_a$ 
and  $\mcalThat_a$ by 
putting for all $\bfr=(\bfr_{\bfv},\bfr_p)^{T}\in Z$
and all $\bfs=(\bfs_{\bfv},\bfs_p)^{T}\in \Zhat$
\begin{equation}
   \label{Atrm.8a}
\big(\mcalT_a \bfr\big)(t)\,=\, 
\left( {-\bfr_{\bfv}(2T-t) \atop \bfr_p(2T-t)}\right),\quad 
\big(\mcalThat_a \bfs\big)(t)\,=\, 
\left( {-\bfs_{\bfv}(2T-t) \atop \bfs_p(2T-t)}\right)
\end{equation}
\end{definition}

The following lemma is easy to verify.
\begin{lemma}
    \label{Lemma.Atrm.1}
We have the following commutations
\begin{equation} \label{Atrm.9}
\begin{array}{rcl}
M \mcalS_a=\mcalT_a M,  & \qquad & \mcalS_a Q= Q \mcalT_a, \\ 
\Mhat \mcalS_a=\mcalThat_a \Mhat, 
 & \qquad & \mcalS_a \Qhat= \Qhat \mcalThat_a.
\end{array}
\end{equation}
\end{lemma}

\begin{theorem}
    \label{Theorem.Atrm.2}
For $(\phi,\psi)^T\in U$ we have
\begin{equation} 
   \label{Atrm.10}
\Gamma^{-1} \mcalS_a F \mcalS_a \Gamma \left( {\phi \atop \psi}\right)
\,=\,F^{*} \left( {\phi \atop \psi}\right)
\end{equation}
\end{theorem}

\noindent {\sl Proof:}\quad For the proof it is convenient to
 make the following definition. 
\begin{definition}
 {\bf Acoustic time-reversal experiment:}
 For given $(\phi,\psi)^T\in U$ define $\bfq_{\bfv}(\bfx,t)=\rho\phi$ and 
$\bfq_{p}(\bfx,t)=\kappa \psi$   and
perform the following physical experiment
\begin{eqnarray} 
   \label{Atrm.11}
\rho\frac{\partial}{\partial s}\bfv_{tr}(\bfx,s) \,+\, \ngrad  
p_{tr}(\bfx,s)  &=& -\bfq_{\bfv}(\bfx,2T-s) \\
   \label{Atrm.12}
\kappa\frac{\partial}{\partial s}p_{tr}(\bfx,s)\,+\, \ndiv  
\bfv_{tr}(\bfx,s) &=& \bfq_{p}(\bfx,2T-s)
\end{eqnarray} 
\begin{equation}
   \label{Atrm.13} 
\bfv_{tr}(\bfx,T)\,=\, 0,\qquad p_{tr}(\bfx,T)\,=\,0\quad \mbox{in}
\quad \Omega,
\end{equation} 
with $s\in [T,2T]$ and zero boundary conditions.
 Doing this experiment means to process the 
data in the following way: Time-reverse all data $(\phi,\psi)^T$ according to 
$t \rightarrow 2T-s$, $t\in [0,T]$, 
and, in addition, reverse the directions of 
the velocities $\phi \rightarrow -\phi$. 
 We call this experiment
the 'acoustic time-reversal experiment'.
Notice that the time is running forward in this experiment.
\end{definition}

The solution $(\bfv_{tr},p_{tr})^T$ 
of this experiment can obviously be represented by 
\begin{equation}
   \label{Atrm.14}
\left( {\bfv_{tr} \atop p_{tr}} \right)\,=\,
F\mcalS_a \Gamma \left( {\phi \atop \psi}\right).
\end{equation}
In order to show that the so defined  time-reversal
experiment is correctly modelled by the adjoint system
derived above, we make the following change in variables: 
\begin{equation}
   \label{Atrm.15} 
\tau\,=\, 2T-s,\quad \hat{\bfv}_{tr}\,=\,-\bfv_{tr},\quad  
\hat{p}_{tr}\,=\, p_{tr}, 
\end{equation} 
which just corresponds to the application of the operator
$\mcalS_a$ to $(\bfv_{tr},p_{tr})^T$.
We have $\tau\in[0,T]$. In these variables the time-reversal system 
(\ref{Atrm.11})-(\ref{Atrm.13}) gets the form 
\begin{eqnarray} 
   \label{Atrm.16}
-\rho\frac{\partial}{\partial \tau} \hat{\bfv}_{tr}(\bfx,\tau)   
\,-\, \ngrad   
\hat{p}_{tr}(\bfx,\tau) &=& \bfq_{\bfv}(\bfx,\tau) \\ 
   \label{Atrm.17}
-\kappa\frac{\partial}{\partial \tau}\hat{p}_{tr}(\bfx,\tau)\,-\, \ndiv   
\hat{\bfv}_{tr}(\bfx,\tau) &=& \bfq_p(\bfx,\tau) 
\end{eqnarray} 
\begin{equation}
   \label{Atrm.18} 
\hat{\bfv}_{tr}(\bfx,T)\,=\, 0,\qquad \hat{p}_{tr}(\bfx,T)\,=\,0 \quad
\mbox{in} \quad \Omega.
\end{equation}
Taking into account the definition of $\bfq_{\bfv}$ and $\bfq_p$, we
see that 
\begin{displaymath}
\left({\hat{\bfv}_{tr} \atop \hat{p}_{tr}}\right)\,=\,
\left({\bfv_{a} \atop p_{a}}\right)
\end{displaymath}
where $\bfv_a$ and $p_a$ solve the adjoint system 
(\ref{Atrm.4})--(\ref{Atrm.6}). Therefore, according to 
theorem \ref{Theorem.Atrm.1}:
\begin{displaymath}
\left({\hat{\bfv}_{tr} \atop \hat{p}_{tr}}\right)\,=\,\Gamma F^{*}
 \left({\phi \atop \psi}\right).
\end{displaymath}
Since we have with (\ref{Atrm.14}), (\ref{Atrm.15})  also 
\begin{displaymath}
\left({\hat{\bfv}_{tr} \atop \hat{p}_{tr}}\right)\,=\,
\mcalS_a F\mcalS_a \Gamma \left( {\phi \atop \psi}\right),
\end{displaymath}
the theorem is proven.  \hfill $\Box$\\

\section{The electromagnetic time-reversal mirror}\label{ETRM.Sec} 

In the following we want to define an operator $\mcalS_e$ 
 such that $F^{*}=\Gamma^{-1}\mcalS_e F \mcalS_e\Gamma$ holds.
 We will call this
operator $\mcalS_e$ the '{\bf electromagnetic time-reversal operator}'.
We will also define the '{\bf electromagnetic time-reversal mirrors}'
 $\mcalT_e$ and  $\mcalThat_e$,
which act on the signals instead of the sources or fields.

We consider Maxwell's equations 
\begin{eqnarray} 
   \label{Etrm.1}
 \epsilon \frac{\partial \bfE_{f}}{\partial t}\,-\, \ncurl \,\bfH_{f}  
\,+\, \sigma\bfE_{f}&=& \bfq_E  \\ 
   \label{Etrm.2}
 \mu \frac{\partial \bfH_{f}}{\partial t}\,+\, \ncurl \,\bfE_{f}  
\quad \qquad &=& \bfq_H 
\end{eqnarray} 
\begin{equation}
   \label{Etrm.3} 
\bfE_{f}(\bfx,0)=0,\qquad \bfH_{f}(\bfx,0)=0  
\end{equation} 
with  $t\in [0,T]$ and zero boundary conditions.
We want to calculate the action of the {\sl adjoint operator}
$F^{*}$
on a vector $(\phi,\psi)^T\in U$. 

\begin{theorem}\label{Theorem.Etrm.1}
Let $(\phi,\psi)^T\in U$ and let $(\bfE_a,\bfH_a)^T$
 be the solution of the
adjoint system
\begin{eqnarray}
   \label{Etrm.4} 
-\epsilon\frac{\partial \bfE_a}{\partial t}\,+\, \ncurl  
\bfH_a(\bfx,t)
\,+\, \sigma\bfE_{a} &=& \epsilon(\bfx)    \phi(\bfx,t)   \\
   \label{Etrm.5}
-\mu\frac{\partial \bfH_a}{\partial t}\,-\, \ncurl 
\bfE_a(\bfx,t)
\qquad \quad &=& \mu(\bfx)    \psi(\bfx,t)  
\end{eqnarray} 
\begin{equation}
   \label{Etrm.6} 
\bfE_{a}(\bfx,T)\,=\, 0,\qquad \bfH_{a}(\bfx,T)\,=\,0\quad
\mbox{{\rm in}}\quad \Omega
\end{equation} 
with  $t\in [0,T]$ and zero boundary conditions. 
Then we have
\begin{equation}
   \label{Etrm.7} 
F^{*}\left( {\phi \atop \psi}\right)\,=\,\left( {\epsilon^{-1} 
 \bfE_a(\bfx,t) \atop \mu^{-1}\bfH_a(\bfx,t)}\right).
\end{equation} 
\end{theorem}
\noindent {\sl Proof:}\quad This theorem is just an application of 
theorem \ref{Theorem.Adj.2} to the electromagnetic symmetric hyperbolic
system. Again, we will give a direct
proof in appendix D  as well.

\begin{definition}
We define the {\bf electromagnetic time-reversal operator} $\mcalS_e$ by 
putting for all $\bfq=(\bfq_{E},\bfq_{H})^{T}\in U$
(and similarly for all $\bfu=(\bfE,\bfH)^{T}\in U$)
\begin{equation}
   \label{Etrm.8}
\big(\mcalS_e \bfq\big)(\bfx,t)\,=\,
 \left({-\bfq_{E}(\bfx,2T-t) \atop \bfq_H(\bfx,2T-t)}\right)
\end{equation}
\end{definition}
\begin{definition}
We define the {\bf electromagnetic time-reversal mirrors} $\mcalT_e$
and   $\mcalThat_e$   by 
putting for all $\bfr=(\bfr_{E},\bfr_{H})^{T}\in Z$
and for all  $\bfs=(\bfs_{E},\bfs_{H})^{T}\in \Zhat$
\begin{equation}
   \label{Etrm.8a}
\big(\mcalT_e \bfr\big)(t)\,=\,
 \left({-\bfr_{E}(2T-t) \atop \bfr_H(2T-t)}\right), \quad
\big(\mcalThat_e \bfs\big)(t)\,=\,
 \left({-\bfs_{E}(2T-t) \atop \bfs_H(2T-t)}\right) \quad
\end{equation}
\end{definition}

The following lemma is easy to verify.
\begin{lemma} \label{Lemma.Etrm.1}
We have the following commutations
\begin{equation}
       \label{Etrm.9}
\begin{array}{rcl}
M \mcalS_e=\mcalT_e M,  & \qquad & \mcalS_e Q= Q \mcalT_e, \\ 
\Mhat \mcalS_e=\mcalThat_e \Mhat, 
 & \qquad & \mcalS_e \Qhat= \Qhat \mcalThat_e.
\end{array}
\end{equation}
\end{lemma}

\begin{theorem} \label{Theorem.Etrm.2}
For $(\phi,\psi)^T\in U$ we have
\begin{equation}
   \label{Etrm.10}
\Gamma^{-1} \mcalS_e F \mcalS_e \Gamma \left( {\phi \atop \psi}\right)
\,=\,F^{*} \left( {\phi \atop \psi}\right)
\end{equation}
\end{theorem}

\noindent {\sl Proof:}\quad For the proof it is convenient to
 make the following definition. 
\begin{definition}
 {\bf Electromagnetic time-reversal experiment:}
 For a given vector $(\phi,\psi)^T\in U$ 
define $\bfq_{E}(\bfx,t)=\epsilon\phi$ and
$\bfq_{H}(\bfx,t)=\mu\psi$ and   perform 
the physical experiment 
\begin{eqnarray}
   \label{Etrm.11} 
\epsilon \frac{\partial}{\partial s}\bfE_{tr}(\bfx,s)\,-\,\ncurl 
\bfH_{tr}(\bfx,s)
\,+\, \sigma \bfE_{tr}(\bfx,s) &=& -\bfq_{E}(\bfx,2T-s)  \\
   \label{Etrm.12}
\mu\frac{\partial}{\partial s}\bfH_{tr}(\bfx,s)\,+\, \ncurl
\bfE_{tr}(\bfx,s)
\qquad \qquad &=& \bfq_{H}(\bfx,2T-s) 
\end{eqnarray} 
\begin{equation}
   \label{Etrm.13} 
\bfE_{tr}(\bfx,T)\,=\, 0,\qquad \bfH_{tr}(\bfx,T)\,=\,0
\quad \mbox{{\rm in}}\quad \Omega,
\end{equation} 
with $s\in [T,2T]$ and zero boundary conditions. 
Doing this experiment means to process the 
data in the following way: Time-reverse all data according to 
$t \rightarrow 2T-s$, $t\in [0,T]$,
and, in addition, reverse the directions of 
the electric field component by
$\phi \rightarrow -\phi$.
We call this experiment
the 'electromagnetic time-reversal experiment'.
Notice that the time is running forward in this experiment.
\end{definition}

The solution $(\bfE_{tr},\bfH_{tr})^T$ 
of this experiment can obviously be represented by 
\begin{equation}
    \label{Etrm.14}
\left( {\bfE_{tr} \atop \bfH_{tr}} \right)\,=\,
F\mcalS_e \Gamma \left( {\phi \atop \psi}\right).
\end{equation}
In order to show that the so defined  time-reversal
experiment is correctly modelled by the adjoint system
derived above, we make the following change in variables: 
\begin{equation}
   \label{Etrm.15} 
\tau\,=\, 2T-s,\quad \hat{\bfE}_{tr}\,=\,-\bfE_{tr},\quad  
\hat{\bfH}_{tr}\,=\, \bfH_{tr}, 
\end{equation} 
which just corresponds to the application of the operator $\mcalS_e$ to
$(\bfE_{tr},\bfH_{tr})^T$.
We have $\tau\in[0,T]$. In these variables the time-reversal system 
(\ref{Etrm.11})--(\ref{Etrm.13}) gets the form 
\begin{eqnarray} 
   \label{Etrm.16}
-\epsilon\frac{\partial}{\partial \tau} \hat{\bfE}_{tr}(\bfx,\tau)   
\,+\, \ncurl   
\hat{\bfH}_{tr}(\bfx,\tau) 
\,+\, \sigma \hat{\bfE}_{tr}(\bfx,s)
&=& \bfq_{\bfE}(\bfx,\tau), \\
   \label{Etrm.17}
-\mu\frac{\partial}{\partial \tau}\hat{\bfH}_{tr}(\bfx,\tau)\,-\,\ncurl   
\hat{\bfE}_{tr}(\bfx,\tau)
\qquad\qquad &=& \bfq_{\bfH}(\bfx,\tau)
\end{eqnarray} 
\begin{equation}
   \label{Etrm.18}
\hat{\bfE}_{tr}(\bfx,T)\,=\, 0,\qquad \hat{\bfH}_{tr}(\bfx,T)\,=\,0
\quad\mbox{in}\quad\Omega.
\end{equation}
Taking into account the definition of $\bfq_{E}$ and $\bfq_H$, we
see that 
\begin{displaymath}
\left({\hat{\bfE}_{tr} \atop \hat{\bfH}_{tr}}\right)\,=\,
\left({\bfE_{a} \atop \bfH_{a}}\right)
\end{displaymath}
where $\bfE_e$ and $\bfH_e$ solve the adjoint system  
(\ref{Etrm.4})--(\ref{Etrm.6}). Therefore, according to 
theorem \ref{Theorem.Etrm.1}:
\begin{displaymath}
\left({\hat{\bfE}_{tr} \atop \hat{\bfH}_{tr}}\right)\,=\,\Gamma F^{*}
 \left({\phi \atop \psi}\right).
\end{displaymath}
Since we have with (\ref{Etrm.14}), (\ref{Etrm.15})  also 
\begin{displaymath}
\left({\hat{\bfE}_{tr} \atop \hat{\bfH}_{tr}}\right)\,=\,
\mcalS_e F\mcalS_e \Gamma \left( {\phi \atop \psi}\right),
\end{displaymath}
the theorem is proven.  \hfill $\Box$\\

\begin{remark}
{\rm
For electromagnetic waves, there is formally an alternative way to define
the {\sl electromagnetic time-reversal operator}, namely 
putting for all $\bfq=(\bfq_{E},\bfq_{H})^{T}\in U$
\begin{displaymath}
\big(\mcalS_e \bfq\big)(\bfx,t)\,=\, 
\left({\bfq_{E}(\bfx,2T-t) \atop -\bfq_H(\bfx,2T-t)}\right),
\end{displaymath}
accompanied by the analogous definitions for the {\sl electromagnetic
time-reversal mirrors}.
With these alternative definitions, theorem \ref{Theorem.Etrm.2} holds
true as well, with only very few changes in the proof. 
 Which form to use depends mainly on the
preferred form for modelling applied antenna signals in the
given antenna system. The first formulation directly works with
applied electric 
 currents, whereas the second form is useful for example for magnetic
dipole sources.
}
\end{remark}

\section{Time-reversal and the adjoint operator $A^{*}$}\label{TRA.Sec}
Define $\mcalS=\mcalS_a$, $\mcalT=\mcalT_a$, $\mcalThat=\mcalThat_a$ 
for the acoustic case, and 
 $\mcalS=\mcalS_e$, $\mcalT=\mcalT_e$, $\mcalThat=\mcalThat_e$
for the electromagnetic case. We call $\mcalS$ the {\bf time-reversal operator}
and $\mcalT$, $\mcalThat$ the  {\bf time-reversal mirrors}.
We combine the results of lemma \ref{Lemma.Atrm.1} and lemma \ref{Lemma.Etrm.1}
into the following lemma. 
\begin{lemma} 
    \label{Lemma.Tra.1}
We have the following commutations 
\begin{equation}
   \label{Tra.1}
\begin{array}{rcl}
M \mcalS=\mcalT M,  & \qquad & \mcalS Q= Q \mcalT, \\ 
\Mhat \mcalS=\mcalThat \Mhat, 
 & \qquad & \mcalS \Qhat= \Qhat \mcalThat.
\end{array}
\end{equation}
\end{lemma}
Moreover, combining theorem \ref{Theorem.Atrm.2} and theorem 
\ref{Theorem.Etrm.2} we get
\begin{theorem} 
   \label{Theorem.Tra.1}
For $(\phi,\psi)^T\in U$ we have
\begin{equation}
   \label{Tra.2}
\Gamma^{-1} \mcalS F \mcalS \Gamma \left( {\phi \atop \psi}\right)
\,=\,F^{*} \left( {\phi \atop \psi}\right).
\end{equation}
\end{theorem} 
With this, we can prove the following theorem which provides the
fundamental link between time-reversal and inverse problems.
\begin{theorem}
   \label{Theorem.Tra.2}
We have
\begin{equation}
   \label{Tra.3}
A^{*}\,=\,\mcalT B \mcalThat.
\end{equation}
\end{theorem}
\noindent {\sl Proof:}\quad 
Recall that the adjoint operator $A^{*}$ can be decomposed
as $A^{*}=Q^{*}F^{*}\Mhat^{*}$. 
With theorem   \ref{Theorem.Adj.1}, theorem \ref{Theorem.Tra.1},
and lemma \ref{Lemma.Tra.1},  it follows therefore 
that 
\begin{eqnarray*}
A^{*}&=&M\Gamma\Gamma^{-1}\mcalS F \mcalS \Gamma\Gamma^{-1}\Qhat \\
&=& M \mcalS F \mcalS \Qhat\\
&=& \mcalT M F \Qhat \mcalThat \\
&=& \mcalT B \mcalThat,
\end{eqnarray*}
which proves the theorem. \hfill $\Box$

\begin{remark}
{\rm The  above  theorem  provides a direct link between the
adjoint operator $A^{*}$, which plays a central role in the
theory of inverse problems,  and a physical 
experiment modelled by $B$.
The expression $\mcalT B \mcalThat$ defines a 'time-reversal experiment'.
We will demonstrate in the following sections
 how we can make use of this relationship in order
to solve the inverse problem of communication by a series of physical 
time-reversal experiments.
} 
\end{remark}
\begin{remark}
{\rm
We only mention here that the above results hold as well for  
{\bf elastic waves} with a suitable
definition of the {\bf elastic time-reversal mirrors}. We leave
out the details for brevity.  
}
\end{remark}

\section{Iterative time-reversal for the gradient method}\label{GMC.Sec}

\subsection{The basic version}\label{GMC.Sec.1}

The results achieved above give rise to the
 following experimental procedure for applying the gradient
method  (\ref{Mat.4}) to
the inverse problem of communication  as
formulated in sections  \ref{Opt.Sec} and \ref{Mat.Sec}.
First, the  pilot signal $\tilde{\bfs}(t)$ is defined
by user $\mcalU_1$ as described in (\ref{Mat.1}). Moreover, we
assume that the first guess $\bfr^{(0)}(t)$
at the base station is chosen to be zero. Then, using 
theorem \ref{Theorem.Tra.2},
we can write  the gradient method 
 (\ref{Mat.4})  in the equivalent form 
\begin{equation}
   \label{Gmc.1}
\begin{array}{rcl}
\bfr^{(0)}&=&0\\
\bfr^{(n+1)}&=& \bfr^{(n)}\,+\,\beta^{(n)} \mcalT B \mcalThat
 (\tilde{\bfs}- A \bfr^{(n)}),
\end{array}
\end{equation}
or, expanding it,
\begin{equation}
  \label{Gmc.2}
\begin{array}{rcl}
\bfr^{(0)}&=&0\\
\bfs^{(n)}&=&  \tilde{\bfs}- A \bfr^{(n)}\\
\bfr^{(n+1)}&=& \bfr^{(n)}\,+\,\beta^{(n)} \mcalT B  \mcalThat \bfs^{(n)}
\end{array}
\end{equation}
In a more detailed form, we arrive at the following {\bf experimental
procedure for implementing the gradient method},
 where  we fix in this description $\beta^{(n)}=1$
for all iteration numbers $n$ for simplicity.

\begin{enumerate}
\item The user $\mcalU_1$ chooses a pilot signal $\alpha(t)$ which
he wants to use for communicating with the base station. The
objective signal at all users 
is then $\tilde{\bfs}(t)=(\alpha(t),0,\ldots,0)^T$. The initial
guess $\bfr^{(0)}(t)$ at the base station
is defined to be zero, such that $\bfs^{(0)}=\stilde$.
\item The user $\mcalU_1$  initiates
the communication by sending the time-reversed pilot signal 
into the environment. This signal is 
$\mcalThat\bfs^{(0)}(t)$. All other
users are quiet.  
\item The base station receives the pilot signal as $B\mcalThat\bfs^{(0)}(t)$.
It time-reverses this signal and sends this time-reversed form,
namely $\bfr^{(1)}(t)= \mcalT B \mcalThat \bfs^{(0)}(t)$, 
back into the medium.  
\item The new signal arrives at all users as $A\bfr^{(1)}(t)$.
All users compare the received signals with their components of the
objective signal  $\tilde{\bfs}(t)$. They take the difference 
$\bfs^{(1)}(t)= \tilde{\bfs}(t)-  A\bfr^{(1)}(t) $,  and
time-reverse it.  They send this new signal 
$\mcalThat \bfs^{(1)}(t)$  back into the medium.
\item The base station receives this new signal, time-reverses it,
adds it to the previous signal $\bfr^{(1)}(t)$,
and sends the sum back into the medium as $\bfr^{(2)}(t)$.
\item This iteration is continued until all users are satisfied
with the match between the received signal $ A\bfr^{(n)}(t)$
and the objective signal $\tilde{\bfs}(t)$ at some iteration 
number $n$. Alternatively, a
fixed iteration number $n$ can be specified a-priori for stopping
the iteration.  
\end{enumerate}
Needless to say that, in practical implementations, the laboratory time needs
to be reset to zero after each time-reversal step. 

\begin{remark}
{\rm 
The experimental procedure which is described above is practically
equivalent to the experimental procedure which was suggested and experimentally
verified in \cite{Montaldo04a,Montaldo04b}. Therefore, our
basic scheme provides an alternative derivation and interpretation of this
experimental procedure.
}
\end{remark}

\begin{remark}
{\rm
We mention that several refinements of this scheme are 
possible and straightforward. For example, a weighted inner
product can be introduced for the user signal space $\Zhat$
 which puts different preferences
 on the satisfaction of the user objectives during
the iterative optimization process. For example, if the 
'importance' of suppressing interferences with other users 
is valued higher than to get an optimal signal quality at 
the specified user $\mcalU_1$, 
a higher weight can be put into the inner product at those
users which did not start the communication process. 
A user who does not care about these interferences, simply puts
 a very small weight into his component of the
inner product of $\Zhat$.
}
\end{remark}

\begin{remark}
{\rm
Notice that there is no mechanism directly built into this procedure
which prevents the energy emitted by the base antenna 
to increase more than the communication system
can support.
For example, if the subspace of signals 
\begin{displaymath}
Z_0\,:=\,\{ \bfr(t)\,:\, A  \bfr = 0 \}
\end{displaymath}
is not empty, then it might happen that during the iteration 
described above (e.g. due to  noise)
 an increasing amount of energy is put
into signals emitted by the base station which are in this subspace 
and which all produce zero contributions to the measurements
at all users. 
More generally, elements of the subspace of signals 
\begin{displaymath}
Z_{\varepsilon}\,:=\,\{ \bfr(t)\,:\, 
 \| A  \bfr\|_{\Zhat}< \varepsilon  \| \bfr \|_{Z} \},
\end{displaymath}
for a very small threshold $0<\epsilon<<1$,
might cause problems during the iteration if the pilot 
signal $\tilde{\bfs}(t)$ chosen by the user has contributions 
in the subspace $A Z_{\varepsilon}$ (i.e. in the space of all
$\bfs=A\bfr$ with $\bfr\in Z_{\varepsilon}$). This is so because 
in the effort of decreasing the mismatch between 
$A\bfr$ and  $\tilde{\bfs}(t)$, the base antenna might
need to put signals with high energy into the system
in order to get only small improvements in the signal match 
at the user side. 
Since the environment (and therefore the operator
$A$) is unknown a-priori, it is difficult to avoid
the existence of such contributions in the pilot signal.
}
\end{remark}

One possible way to prevent the energy
emitted by the base station to increase 
artificially would be to project the signals $\bfr^{(n)}(t)$
onto the orthogonal complements of the subspaces
$Z_0$ or $Z_{\varepsilon}$ (if they are known or can be constructed
by some means)
 prior to their emission. 
Alternatively,  the iteration can be stopped at an early stage 
before these unwanted contributions start to build up. (This in fact
has been suggested in \cite{Montaldo04a,Montaldo04b}).

In the following subsection we introduce an alternative way of 
ensuring that the energy emitted by the base station stays
reasonably bounded in the effort of fitting the pilot signal
at the users.

\subsection{The regularized version}\label{GMC.Sec.2}

Consider the regularized problem
\begin{equation}
    \label{Gmc.3}
\bfr_{LSr}\,=\, \mbox{Min}_{\bfr}\big(
\| A\bfr - \stilde\|_{\Zhat}^2
\,+\, \lambda \|\bfr\|_{Z}^2\big)
\end{equation}
with some suitably chosen regularization parameter $\lambda>0$.
In this problem formulation a trade-off is sought between
a signal fit at the user side and a minimized energy emission
at the base station. The trade-off parameter is the regularization
parameter $\lambda$.
Instead of (\ref{Mat.3}) we need to consider now 
\begin{equation} 
    \label{Gmc.4}
\tilde{\mcalJ}(\bfr)\,=\,\frac{1}{2}\Big\langle A\bfr - \tilde{\bfs} 
\,,\, A\bfr - \tilde{\bfs}\Big\rangle_{\Zhat}\,+\, 
\frac{\lambda}{2}\Big\langle \bfr  \,,\, \bfr \Big\rangle_{Z}
\end{equation}
The negative gradient direction is now given by 
$-A^{*}(A\bfr-\tilde{\bfs})-\lambda \bfr$, such that the 
regularized iteration reads:
\begin{equation}
  \label{Gmc.5}
\begin{array}{rcl}
\bfr^{(0)}&=&0\\
\bfr^{(n+1)}&=& \bfr^{(n)}\,+\,\beta^{(n)}
\Big(\mcalT B \mcalThat (\tilde{\bfs} - A\bfr^{(n)})-\lambda \bfr^{(n)}\Big),
\end{array}
\end{equation}
where we have replaced $A^{*}$ by $\mcalT B \mcalThat $. The time-reversal 
iteration can be expanded into the following practical scheme
\begin{equation}
  \label{Gmc.6}
\begin{array}{rcl}
\bfr^{(0)}&=&0\\
\bfs^{(n)}&=&  \tilde{\bfs}- A \bfr^{(n)}\\
\bfr^{(n+1)}&=& \bfr^{(n)}\,+\,\beta^{(n)} \mcalT B  \mcalThat \bfs^{(n)}
\,-\,\beta^{(n)} \lambda \bfr^{(n)}.
\end{array}
\end{equation}
Comparing with (\ref{Gmc.2}),  we see that the adaptations which 
need to be applied
in the practical implementation for stabilizing the 
basic algorithm can easily be done.

\section{The minimum norm solution approach}\label{MNS.Sec}

\subsection{The basic version}\label{MNS.Sec.1}

In this section we want to propose an alternative scheme
for solving the inverse problem of communication. As mentioned above,
a major drawback of the basic approach (\ref{Mat.2}) is that the
energy emitted by the base station is not limited explicitly
when solving the optimization problem. The regularized version
presented above alleviates this problem. However, we want
to mention here that, under certain assumptions, 
there is an alternative scheme which can be employed instead
and which has an energy constraint directly built in. 
Under the formal assumption that there exists at least one (and presumably
more than one) solution
of the inverse problem at hand (i.e. the 'formally
underdetermined case'), we can look for the 
{\bf minimum norm solution}
\begin{equation}
   \label{Mns.1}
\mbox{Min}_{\bfr} \|\bfr\|_{Z}\quad
\mbox{subject to}\quad
A\bfr\,=\,\tilde{\bfs}.
\end{equation}
In Hilbert spaces this solution has an explicit form. It is
\begin{equation}
   \label{Mns.2}
\bfr_{MN}\,=\,\Astar (A\Astar)^{-1}\tilde{\bfs}.
\end{equation}
Here, the operator $(A\Astar)^{-1}$ acts as a filter on 
the pilot signal $\stilde$. Instead of sending the pilot
signal to the base station, the users send the filtered
version of it. Certainly, a method must be found in order
to apply the filter $(A\Astar)^{-1}$ to the pilot signal. 
One possibility of doing so would be to try to determine the
operator $A\Astar$ explicitly by a series of time-reversal
experiments on some set of basis functions of $\Zhat$, and then invert this
operator numerically. However, this might not be practical 
in many situations. (It certainly would be slow
and it would involve a significant
amount of  signal-processing, 
which we want to avoid here.) 
Therefore, we propose an alternative procedure. First, we
notice that there is no need to determine the whole
operator  $(A\Astar)^{-1}$, but that we only have to apply it to 
one specific signal, namely $\stilde$.  
Let us introduce the short notation
\begin{displaymath}
C\,=\,A\Astar.
\end{displaymath}
In this notation, we are
looking for a signal $\scirc\in \Zhat$ such that $C \scirc=\stilde$.
We propose to solve this equation in the least squares sense:
\begin{equation}
   \label{Mns.3}
\scirc\,=\,\mbox{Min}_{\bfs} \|C\bfs -\stilde\|_{\Zhat}^2.
\end{equation}
Moreover, as a suitable method for practically finding this
solution, we want to use the {\bf gradient method}. Starting
with the initial guess $\bfs_0=0$, the gradient method reads
\begin{equation}
  \label{Mns.4}
\begin{array}{rcl}
\bfs^{(0)}&=&0\\
\bfs^{(n+1)}&=& \bfs^{(n)}\,-\,\beta^{(n)} C^{*}
 (C \bfs^{(n)}- \stilde),
\end{array}
\end{equation}
where $C^{*}$ is the adjoint operator to $C$ and $\beta^{(n)}$ is
again some step-size.
Expanding this expression, and taking into account $C^{*}=C$ and
$\Astar=\mcalT B \mcalThat$, we arrive at
\begin{equation}
  \label{Mns.5}
\begin{array}{rcl}
\bfs^{(0)}&=&0\\
\bfs^{(n+1)}&=& \bfs^{(n)}\,+\,\beta^{(n)} A\mcalT B \mcalThat 
 \big( \stilde - A\mcalT B \mcalThat \bfs^{(n)} \big).
\end{array}
\end{equation}
In the practical implementation, we arrive at the following
 iterative scheme:
\begin{equation}
\fl
  \label{Mns.6}
\begin{array}{rcl}
\mbox{Initialize gradient iteration:} & & \bfs^{(0)}=  0\\
\mbox{for}\,\,n=0,1,2,\ldots \,\, \mbox{do:}  & & \\
\bfr^{(n+\frac{1}{2})} &=& \mcalT B \mcalThat \bfs^{(n)}\\
\bfs^{(n+\frac{1}{2})}&=&  \stilde - A \bfr^{(n+\frac{1}{2})}\\
\bfr^{(n+1)}&=& \mcalT B \mcalThat \bfs^{(n+\frac{1}{2})}\\
\bfs^{(n+1)}&=& \bfs^{(n)}+ \beta^{(n)} A \bfr^{(n+1)}\\
\mbox{Terminate iteration at step } \nmax \mbox{:} & &  \scirc = \bfs^{(\nmax)} \\
\mbox{Final result:} & &  \bfr_{MN}=\mcalT B \mcalThat \scirc.
\end{array}
\end{equation}
This iteration can be implemented by a series of time-reversal
experiments, without the need of heavy signal-processing. 
The final
step of the above algorithm amounts to applying 
$A^{*}$ to the result of the gradient
iteration for calculating $\scirc=(A\Astar)^{-1}\stilde$,
 which  yields then $\bfr_{MN}$. This will then be the signal
to be applied by the base station during the communication process with
the user $\mcalU_1$.

\subsection{The regularized version}\label{MNS.Sec.2}

In some situations it might be expected that the operator $C$ is
ill-conditioned, such that its inversion might cause instabilities,
in particular when noisy signals are involved. For those situations,
a regularized form of the minimum norm solution is available, namely
\begin{equation}
   \label{Mns.7}
\bfr_{MNr}\,=\,\Astar (A\Astar+\lambda I_{\Zhat})^{-1}\tilde{\bfs}
\end{equation}
where $I_{\Zhat}$ denotes the identity operator in $\Zhat$ and $\lambda>0$
is some suitably chosen regularization parameter. The necessary
adjustments in the gradient iteration for applying 
$(A\Astar+\lambda I_{\Zhat})^{-1}$
to $ \tilde{\bfs}$
are easily done. We only
mention here the resulting procedure for the implementation of this 
gradient method by 
a series of time-reversal experiments:
\begin{equation}
\fl
  \label{Mns.8}
\begin{array}{rcl}
\mbox{Initialize gradient iteration:} & & \bfs^{(0)}=  0\\
\mbox{for}\,\,n=0,1,2,\ldots \,\, \mbox{do:}  & & \\
\bfr^{(n+\frac{1}{2})} &=& \mcalT B \mcalThat \bfs^{(n)}\\
\bfs^{(n+\frac{1}{2})}&=&  \stilde - A \bfr^{(n+\frac{1}{2})}-\lambda \bfs^{(n)}\\
\bfr^{(n+1)}&=& \mcalT B \mcalThat \bfs^{(n+\frac{1}{2})}\\
\bfs^{(n+1)}&=& \bfs^{(n)}+ \beta^{(n)} A \bfr^{(n+1)}+
 \beta^{(n)}\lambda \bfs^{(n+\frac{1}{2} )}  . \\
\mbox{Terminate iteration at step } \nmax \mbox{:} & &  \scirc = \bfs^{(\nmax)} \\
\mbox{Final result:} & &  \bfr_{MNr}=\mcalT B \mcalThat \scirc.
\end{array}
\end{equation}
Again, the last step shown above is a final application 
of $A^{*}$ to the result of the gradient
iteration for calculating $\scirc=(A\Astar+\lambda I_{\Zhat})^{-1}\stilde$,
which yields then $\bfr_{MNr}$. This will then be the signal
to be applied by the base station during the communication process with
the user $\mcalU_1$.

\section{The regularized least squares solution revisited}\label{RLS.Sec}

We have introduced above the  regularized least
squares solution of the inverse problem of communication, namely
\begin{equation}
   \label{Rls.1}
\bfr_{LSr}\,=\, \mbox{Min}_{\bfr}\big(
\| A\bfr - \stilde\|_{\Zhat}^2
\,+\, \lambda \|\bfr\|_{Z}^2\big)
\end{equation}
with $\lambda>0$
being the regularization parameter.
In Hilbert spaces, the solution of (\ref{Rls.1}) has an 
explicit form. It is 
\begin{equation}
   \label{Rls.2}
\bfr_{LSr}\,=\,(\Astar A+\lambda I_{Z})^{-1}\Astar \stilde,
\end{equation}
where $I_{Z}$ is the identity operator in $Z$. 
It is therefore tempting to try to implement also this direct
form as a series of time-reversal experiments and compare
its performance with the gradient method as it was described above.
As our last strategy which we present in this paper, we
 want to show here that such an alternative direct implementation 
of  (\ref{Rls.1})  is in fact possible.  

Notice that in (\ref{Rls.1})
 the filtering operator $ (\Astar A+\lambda I_{Z})^{-1}$ is
applied at the base station, in contrast to the previous case
where the user signal was filtered by the operator 
$ (A \Astar +\lambda I_{\Zhat})^{-1}$. Analogously
to the previous case, we need to find a practical way to 
apply this filter to a signal at the base station. We 
propose again to solve the equation 
\begin{equation}
   \label{Rls.3}
(\Astar A+\lambda I_{Z})\rcirc \,=\, \rtilde
\end{equation}
in the least squares sense, where $\rtilde=\Astar \stilde$.
Defining 
\begin{displaymath}
C\,=\,\Astar A+\lambda I_{Z},
\end{displaymath}
and using $C^{*}=C$ and $\Astar=\mcalT B \mcalThat$,
we arrive at the following {\bf gradient iteration} for solving
problem (\ref{Rls.3}):
\begin{equation}
  \label{Rls.4}
\begin{array}{rcl}
\bfr^{(0)}&=& 0 \\
\bfr^{(n+1)}&=& \bfr^{(n)}\,+\,\beta^{(n)} (\mcalT B \mcalThat A+\lambda I_{Z})
 \Big(\rtilde -  (\mcalT B \mcalThat A+\lambda I_{Z})  \bfr^{(n)}  \Big).
\end{array}
\end{equation}
This gives rise to the
following practical implementation  by a series of
time-reversal experiments:
\begin{equation}
\fl
  \label{Rls.5}
\begin{array}{rcl}
\mbox{User sends pilot signal to base station:} && \rtilde =
 \mcalT B \mcalThat \stilde \\
\mbox{Initialize gradient iteration:} & & \bfr^{(0)}=  0\\
\mbox{for}\,\,n=0,1,2,\ldots \,\, \mbox{do:}  & & \\
\bfs^{(n+\frac{1}{2})} &=& A \bfr^{(n)}\\
\bfr^{(n+\frac{1}{2})}&=& 
 \rtilde - \mcalT B \mcalThat \bfs^{(n+\frac{1}{2})}-\lambda \bfr^{(n)}\\
\bfs^{(n+1)}&=& A \bfr^{(n+\frac{1}{2})}\\
\bfr^{(n+1)}&=& \bfr^{(n)}+ \beta^{(n)} \mcalT B \mcalThat  \bfs^{(n+1)}+ 
\beta^{(n)}\lambda \bfr^{(n+\frac{1}{2} )}  . \\
\mbox{Terminate iteration at step } \nmax \mbox{:}  & &  
\rcirc = \bfr^{(\nmax)} \\
\mbox{Final result:} & &  \bfr_{LSr}= \rcirc.
\end{array}
\end{equation}
$\bfr_{LSr}$ will then be the signal
to be applied by the base station during the communication process with
the user $\mcalU_1$.

\section{Partial and generalized measurements} \label{PM.Sec}

In many practical applications, only partial measurements
of the whole wave-field are available. For example, in 
ocean acoustics often only pressure is measured, whereas
the velocity field is not part of the measurement process. 
Similarly, in wireless communication 
only one or two components of the electric field might be
measured simultaneously, but the remaining electric components
and all magnetic components are missing. We want to demonstrate
in this section that all results presented above are valid
also in this situation of partial measurements, with the
suitable adaptations. 

Mathematically, the measurement operator needs to be adapted
for the situation of partial measurements. Let us concentrate
here on the special situation that only one component $\bfu_{\nu}$
$(\nu\in\{1,2,3,...\})$ of the incoming wave field
$\bfu$ is  measured by the users and 
the base station. All other possible situations will then just
be combinations of this particular case. 
It might also
occur the situation that users can measure a different partial set of 
components than the base station. That case also follows directly
from this canonical situation. 

We introduce the new signal space at the base station 
$Y=(L_2[0,T])^{\anzant}$ and the corresponding
    'signal projection operator'
$P_{\nu}$ by putting
\begin{displaymath}
P_{\nu}\,:\,Z\rightarrow  Y \quad , \quad
P_{\nu}(\bfr_1(t),\ldots,\bfr_{\anzant}(t))^T\,=\, \bfr_{\nu}(t).
\end{displaymath}
We see immediately that its adjoint $P_{\nu}^{*}$ is given by
\begin{displaymath}
P_{\nu}^{*}\,:\,Y\rightarrow  Z \quad , \quad
P_{\nu}^{*}\bfr_{\nu}(t)\,=\,(0,\ldots,0,\bfr_{\nu}(t),0,\ldots,0)^T
\end{displaymath}
where $\bfr_{\nu}(t)$ appears on the right hand side at the $\nu$-th position.
Our new measurement operator $M_{\nu}$, and the new source
operator $Q_{\nu}$, are then defined by
\begin{equation}
  \label{Pm.3} 
\begin{array}{r@{\,=\,}l}
M_{\nu}\,:\,U\rightarrow  Y  , \qquad
M_{\nu}\bfu& P_{\nu}M\bfu\\ 
Q_{\nu}\,:\,Y\rightarrow U, \qquad
Q_{\nu} \bfr_{\nu}  & Q P_{\nu}^{*}\bfr_{\nu}.
\end{array}
\end{equation}
Analogous definitions are done for $\hat{Y}$, $\hat{P}_{\nu}$,
$\Mhat_{\nu}$ and $\Qhat_{\nu}$ at the users.

Obviously, we will have to replace now in the above derivation
of the iterative time-reversal procedure all measurement
operators $M$ by $M_{\nu}$ (and $\Mhat$ by $\Mhat_{\nu}$)
 and all source operators
$Q$ by $Q_{\nu}$ (and $\Qhat$ by $\Qhat_{\nu}$). 
In particular, the new 'communication operators' are now given by 
\begin{equation}
   \label{Pm.4} 
\begin{array}{c@{\quad}c}
A_{\nu}:\,Y\rightarrow \Yhat, & A_{\nu}\bfr_{\nu}
\,=\,\Mhat_{\nu} FQ_{\nu}\bfr_{\nu}, \\
B_{\nu}:\,\Yhat \rightarrow Y, & B_{\nu}\bfs_{\nu}
\,=\,M_{\nu} F \Qhat_{\nu} \bfs_{\nu}.
\end{array}
\end{equation}
 In the following
two theorems we show that the main results of this paper
carry over to these newly defined operators. 

\begin{theorem}
   \label{Theorem.Pm.1}
We have 
\begin{equation}
   \label{Pm.5} 
\begin{array} {c@{\quad}c}
M_{\nu}^{*}\,=\,\Gamma^{-1}Q_{\nu}, &
\Mhat_{\nu}^{*}\,=\,\Gamma^{-1}\Qhat_{\nu},\\
Q_{\nu}^{*}\,=\,M_{\nu} \Gamma , &
\Qhat_{\nu}^{*}\,=\,\Mhat_{\nu} \Gamma.
\end{array}
\end{equation}
\end{theorem}
\noindent {\sl Proof:} \quad The proof is an easy exercise using 
(\ref{Pm.3}) and theorem \ref{Theorem.Adj.1}. \hfill $\Box$

\begin{theorem}
   \label{Theorem.Pm.2}
It is 
\begin{equation}
A_{\nu}^{*}\,=\,\mcalT B_{\nu} \mcalThat.
\end{equation}
\end{theorem}
\noindent {\sl Proof:} \quad The proof is now identical to the proof
of theorem \ref{Theorem.Tra.2}, using theorem \ref{Theorem.Pm.1}
instead of  theorem \ref{Theorem.Adj.1}.  \hfill $\Box$

\begin{remark}
{\rm
In fact, it is easy to verify that all results of this
paper remain
valid for {\sl arbitrarily defined linear measurement
operators} 
\begin{displaymath}
M_{\mcalU}\,:\,U\rightarrow Z_{\mcalU}, \qquad 
M_{\mcalA}\,:\,U\rightarrow Z_{\mcalA}, \qquad 
\end{displaymath}
where $Z_{\mcalU}$ and  $Z_{\mcalA}$ are any meaningful 
signal spaces at the 
users and the base antennas, respectively. The only
requirement is that it is experimentally possible to apply
signals according to the source operators defined by
\begin{eqnarray*}
Q_{\mcalU}\,:\, Z_{\mcalU}\rightarrow U, &\qquad& Q_{\mcalU}=\Gamma M_{\mcalU}^{*}\\
Q_{\mcalA}\,:\, Z_{\mcalA}\rightarrow U,& \qquad&  Q_{\mcalA}=\Gamma M_{\mcalA}^{*}
\end{eqnarray*}
where $M_{\mcalU}^{*}$ and $M_{\mcalA}^{*}$ are the formal adjoint
operators to $M_{\mcalU}$ and $M_{\mcalA}$ with respect
to the chosen signal spaces $Z_{\mcalU}$ and  $Z_{\mcalA}$.
In addition, the measurement and source operators
as defined above are required to
satisfy the commutation relations as stated in 
lemma \ref{Lemma.Tra.1}.
Under these assumptions, we define the {\sl generalized
communication operators} $\tilde{A}$ and $\tilde{B}$ by
\begin{eqnarray*}
\tilde{A}:\,Z_{\mcalA}\rightarrow Z_{\mcalU} , &\quad& 
\tilde{A}\bfr\,=\,M_{\mcalU} F Q_{\mcalA} \bfr, \\
\tilde{B}:\,Z_{\mcalU}\rightarrow Z_{\mcalA}, &\quad&
\tilde{B}\bfs\,=\,M_{\mcalA} F Q_{\mcalU} \bfs.
\end{eqnarray*}
Now the proof to theorem \ref{Theorem.Tra.2} directly carries over to
this generalized situation, such that we have also here
\begin{displaymath}
\tilde{A}^{*}\,=\,\mcalT \tilde{B} \mcalThat.
\end{displaymath}
This yields iterative time-reversal
schemes completely analogous to those presented above.
}
\end{remark}

\section{Summary and future directions}\label{SFR.Sec}

We have derived in this paper a direct link between the
time-reversal technique and the adjoint method for imaging.
Using this relationship, we have constructed several iterative
time-reversal schemes for solving an inverse problem which
arises in ocean acoustic and wireless communication. 
Each of these schemes can be realized physically as
a series of time-reversal experiments, without the use
of heavy signal processing or computations.
 One of the schemes which we have
derived (and which we call the
'basic scheme'), is practically equivalent to a technique introduced
earlier in \cite{Montaldo04a,Montaldo04b} using different
tools. Therefore, we have given an alternative theoretical
derivation of that technique, with a different mathematical
interpretation. The other schemes which we have introduced
are new in this application. They represent either generalizations
of the basic scheme, or alternatives which follow different
objectives. 

Many questions related to these and similar iterative time-reversal
approaches for telecommunication are still open. The experimental
implementation has been investigated so far only for one of these
techniques in  \cite{Montaldo04a,Montaldo04b}, for the situation of 
underwater sound propagation. A thorough experimental
(or, alternatively, numerical) verification of the other schemes
is  necessary for their practical evaluation. 
An interesting and practically important
problem is the derivation of quantitative estimates
 for the expected focusing
quality of each of these schemes, 
for example following the ideas of the work 
performed for a single step in \cite{Blomgren02}. Certainly,
it is expected that these estimates will again strongly depend
on the randomness of the medium, on the geometry and distribution
of users and base antennas, and on technical constraints
as for example partial measurements. Also,
different types of noise in the communication
system need to be taken into account. The performance
in a time-varying environment is another interesting issue of 
practical importance. All schemes presented here can be adapted 
in principle to a dynamic environment by re-adjusting the 
constructed optimal signals periodically. Practical ways of doing so
need to be explored. 

Finally, we  want to mention that speed is an important factor
in the application of communication. 
In the time-reversal experiment, in a certain sense
a physical process is
used for applying a mathematical operator (namely $A^{*}$)
in order to solve a given inverse problem. In a complex environment,
this can be much faster than using a powerful computer, even if
the complex environment would be known to the users and 
the base station.


\ack
The author thanks George Papanicolaou, Hongkai
Zhao and Knut Solna for many exciting and useful discussions on
time-reversal. He thanks Dominique Lesselier for many fruitful 
discussions, and for making the research stay at
Sup\'elec possible. 
He thanks Mathias
Fink, Claire Prada, Gabriel Montaldo, Francois Vignon, 
and the group at Laboratoire Ondes
et Acoustique in Paris 
for pointing out to him their recent work 
on iterative communication schemes, which stimulated
part of this paper.  The author also thanks
Bill Kuperman and the group at Scripps Institution
of Oceanography 
in San Diego for useful discussions. 
Financial support from the Office of Naval Research
(under grant N00014-02-1-0090) and  the CNRS (under grant
 No.~8011618)
is gratefully acknowledged.


\appendix
\setcounter{section}{1}

\section*{Appendix A: Proof of theorem \ref{Theorem.Adj.1}} 

For arbitrary $\bfr\in Z$ and $\bfu\in U$ we have
\begin{eqnarray*}
\left\langle M\bfu, \bfr\right\rangle_{Z}&=&
\sum_{k=1}^{\anzant}\int_{[0,T]}\left\langle
\int_{\Omega}\gamma_k(\bfx)\bfu(\bfx,t)d\bfx\,,\,
\bfr(t)\right\rangle_{N}dt \\
&=&
\int_{[0,T]}\int_{\Omega}\left\langle\Gamma(\bfx)\bfu(\bfx,t)
\,,\,\Gamma^{-1}(\bfx)\sum_{k=1}^{\anzant}\gamma_k(\bfx)
\bfr(t)\right\rangle_{N}d\bfx dt.
\end{eqnarray*}
Therefore,
\begin{displaymath}
\left\langle M\bfu, \bfr\right\rangle_{Z}\,=\,
\left\langle \bfu, M^{*}\bfr\right\rangle_{U}
\end{displaymath}
with
\begin{displaymath}
M^{*}\bfr\,=\,\Gamma^{-1}(\bfx)\sum_{k=1}^{\anzant}\gamma_k(\bfx)
\bfr(t)\,=\,\Gamma^{-1}(\bfx) Q \bfr.
\end{displaymath} 
Doing the analogous calculation for $\Mhat^{*}$ we get
\begin{equation}
  \label{Pt1.1}
M^{*}\,=\,\Gamma^{-1}Q,\quad\quad
\Mhat^{*}\,=\,\Gamma^{-1}\Qhat.
\end{equation}
Taking the adjoint of (\ref{Pt1.1})  we see that
\begin{displaymath}
Q^{*}\,=\,M \Gamma ,\quad\quad
\Qhat^{*}\,=\,\Mhat \Gamma
\end{displaymath}
holds as well. This  proves  the theorem.  \hfill $\Box$\\

\section*{Appendix B: Proof of theorem \ref{Theorem.Adj.2}}

We have the following  version of  {\sl Green's formula}
\begin{equation}  
\fl
   \label{Pt2.1}
\int_{[0,T]}\int_{\Omega} \left\langle
\Gamma(\bfx)\frac{\partial \bfu}{\partial t}\,+\,\sum_{i=1}^{3}
D^i\frac{\partial \bfu}{\partial x_i}\,+\,\Phi(\bfx)\bfu\,,\,\bfz
\right\rangle_N\,d\bfx dt
\end{equation} 
\begin{displaymath}
\fl
+\, \int_{[0,T]}\int_{\Omega} \left\langle
\Gamma(\bfx) \bfu(\bfx,t), \bfv(\bfx,t) \right\rangle_N
\,d\bfx dt
\, + \,
 \int_{[0,T]}\int_{\Omega} \left\langle
 \bfq(\bfx,t)\,,\,\bfz(\bfx,t) 
 \right\rangle_N \,d\bfx dt
\end{displaymath}
\begin{displaymath}
\fl
=\,\int_{[0,T]}\int_{\Omega} \left\langle
\bfu \,,\, - \Gamma(\bfx)\frac{\partial \bfz}{\partial t}
\,-\,\sum_{i=1}^{3}
D^i\frac{\partial \bfz}{\partial x_i}\,+\,\Phi(\bfx)\bfz
\right\rangle_N
\,d\bfx dt
\end{displaymath}
\begin{displaymath}
\fl
+\,\int_{[0,T]}\int_{\Omega}\left\langle \Gamma(\bfx) 
\bfu(\bfx,t),\bfv(\bfx,t)\right\rangle_N\,d\bfx dt
\, +\,
 \int_{[0,T]}\int_{\Omega} \left\langle
 \bfq(\bfx,t)\,,\,\bfz(\bfx,t) 
 \right\rangle_N \,d\bfx dt
\end{displaymath}
\begin{displaymath}
\fl
+\,\int_{\Omega}\left\langle\Gamma(\bfx)\bfu(\bfx,T),\bfz(\bfx,T)
\right\rangle_N\,d\bfx
\,-\,\int_{\Omega}\left\langle\Gamma(\bfx)\bfu(\bfx,0),
\bfz(\bfx,0)\right\rangle_N\,d\bfx
\end{displaymath}
\begin{displaymath}
\fl
+\,\int_{[0,T]}\int_{\partial\Omega}\sum_{i=1}^{3}
\left\langle D^i\bfu,\bfz\right\rangle_N\nu_i(\bfx)
\,d\sigma dt,
\end{displaymath}
where $\bfn(\bfx)=(\nu_1(\bfx),\nu_2(\bfx),\nu_3(\bfx))$ is the
 outward normal at $\partial\Omega$ in the point $\bfx$.
Notice that we have augmented Green's formula
in  (\ref{Pt2.1})  by some terms which appear in identical
form  on the left hand side and on the right hand side.

 We will assume here that the boundary is far away from
the sources and receivers and that no energy enters $\Omega$ 
from the outside, such that
 during the time interval of interest  $[0,T]$ all fields
along this boundary are identically zero. This is 
expressed by the boundary conditions given in 
 (\ref{Sym.4}) and (\ref{Adj.5}).  
Let $\bfu(\bfx,t)$ be a solution of
 (\ref{Sym.1}),(\ref{Sym.2}),(\ref{Sym.4}),
 and $\bfz(\bfx,t)$ a solution
of (\ref{Adj.3})--(\ref{Adj.5}). Then the first 
term on the left hand side of (\ref{Pt2.1}) and the third term on the
right hand side cancel each other because of (\ref{Sym.1}).
The second term on the left hand side and the first term on 
the right hand side cancel each other  because of (\ref{Adj.3}). 
The $(t=T)$-term and the $(t=0)$ term vanish due to 
(\ref{Adj.4}) and (\ref{Sym.2}), respectively, and the 
boundary integral vanishes because of 
(\ref{Sym.4}) and (\ref{Adj.5}). 
The remaining terms (i.e. the third term on the left hand side and the
second term on the right hand side) can be written as
\begin{displaymath}
\left\langle F\bfq,\bfv\right\rangle_{U}\,=\,
\left\langle\bfq, F^{*}\bfv\right\rangle_{U},
\end{displaymath}
with $F^{*}\bfv=\Gamma^{-1}(\bfx) \bfz(\bfx,t)$
 as defined  in (\ref{Adj.6}). \hfill $\Box$

\section*{Appendix C: Direct proof of theorem \ref{Theorem.Atrm.1}}

We prove the lemma by using {\sl Greens formula}:
\begin{equation}
\fl
  \label{Pt3.1}
\int_{0}^{T}\int_{\Omega}\left[ \rho\frac{\partial \bfv_{f}}{\partial
t} \bfv_{a}
\,+\, \ngrad  p_{f} \bfv_{a}
\,+\, \kappa\frac{\partial p_{f}}{\partial t}p_{a}
\,+\,\ndiv \bfv_f p_{a}
\right] \,d\bfx dt
\end{equation}
\begin{displaymath}
\fl
+\,\int_{0}^{T}\int_{\Omega}\left[\rho \bfv_f  \phi\,+\,
  \kappa  p_{f} \psi\right]\,d\bfx dt
+\,\int_{0}^{T}\int_{\Omega}\left[ \bfq_{\bfv}\bfv_a \,+\,
\bfq_p p_a \right]\,d\bfx dt
\end{displaymath}
\begin{displaymath}
\fl
=\, \int_{0}^{T}\int_{\Omega}\left[ 
-\bfv_{f} \rho
\frac{\partial\bfv_{a}}{\partial t}\,-\, p_{f}
\ndiv \bfv_{a}
\,-\, p_{f}\kappa \frac{\partial p_{a}}{\partial t}\,-\,
\bfv_{f} \ngrad  p_{a}
\right]\,d\bfx dt 
\end{displaymath}
\begin{displaymath}
\fl
+\,\int_{0}^{T}\int_{\Omega}\left[\rho \bfv_f  \phi\,+\,
  \kappa  p_{f} \psi\right]\,d\bfx dt
+\,\int_{0}^{T}\int_{\Omega}\left[ \bfq_{\bfv}\bfv_a \,+\,
\bfq_p p_a \right]\,d\bfx dt
\end{displaymath}
\begin{displaymath}
\fl
\,+\,\int_{0}^{T}\int_{\partial\Omega}(\bfv_a\cdot \bfn)p_f\, d\sigma dt
\,+\, 
\int_{0}^{T}\int_{\partial\Omega}(\bfv_f\cdot \bfn)p_a\, d\sigma dt
\end{displaymath}
\begin{displaymath}
\fl
\,+\,\int_{\Omega} \rho\big[(\bfv_{f} \bfv_{a})(\bfx,T)-(\bfv_{f}
 \bfv_{a})(\bfx,0)\big]\,d\bfx 
\,+\,\int_{\Omega} \kappa \big[ (p_{f}p_{a})(\bfx,T)
-(p_{f}p_{a})(\bfx,0)\big]\,d\bfx.
\end{displaymath}
This equation has the form (\ref{Pt2.1}). 
Notice that we have augmented Green's formula
in  (\ref{Pt3.1}), as already shown in (\ref{Pt2.1}),
  by some terms which appear in identical
form  on the left hand side and on the right hand side.

The first term on the left hand side of  equation (\ref{Pt3.1})
and the third term on the right hand side cancel each other
due to (\ref{Atrm.1}),(\ref{Atrm.2}). The second term
on the left hand side and the first term on the right hand side 
cancel each
other because of (\ref{Atrm.4}),(\ref{Atrm.5}). 
The ($t=T$)-terms and the ($t=0$)-terms  
vanish due to (\ref{Atrm.6}), (\ref{Atrm.3}), respectively,  and the boundary 
terms vanish because
of zero boundary conditions. We are left over with the equation
\begin{displaymath}
\int_{0}^{T}\int_{\Omega}\left[\rho \bfv_f  \phi\,+\,
  \kappa  p_{f} \psi\right]\,d\bfx dt
\,=\,\int_{0}^{T}\int_{\Omega}\left[ \bfq_{\bfv}\bfv_a \,+\,
\bfq_p p_a \right]\,d\bfx dt
\end{displaymath}
Defining $F^{*}$ by (\ref{Atrm.7}), this can be written as
\begin{displaymath}
\left\langle F \left( {\bfq_{\bfv} \atop \bfq_p}\right)\,,\, 
\left( {\phi \atop \psi}\right)\right\rangle_U \,=\,
\left\langle \left( { \bfq_{\bfv} \atop  \bfq_p   }\right)\,,\, 
F^{*}\left( {\phi \atop \psi}\right)\right\rangle_U. 
\end{displaymath}
Therefore, $F^{*}$ is in fact the adjoint of $F$, and the 
lemma is proven. \hfill $\Box$\\

\section*{Appendix D: Direct proof of theorem  \ref{Theorem.Etrm.1}}

We prove the lemma by using {\sl Green's formula}:
\begin{displaymath}
\fl
\int_{0}^{T}\int_{\Omega}\bigg[ \epsilon\frac{\partial \bfE_{f}}{\partial
t} \bfE_{a}
\,-\,\ncurl  \bfH_{f} \bfE_{a}
\,+\,\sigma \bfE_f \bfE_a
\,+\, \mu\frac{\partial \bfH_{f}}{\partial t}\bfH_{a}
\,+\,\ncurl \bfE_f \bfH_{a}
\bigg] \,d\bfx dt
\end{displaymath}
\begin{displaymath}
\fl
+\,\int_{0}^{T}\int_{\Omega}\left[ \epsilon\bfE_{f} \phi\,+\,
\mu\bfH_f \psi\right]\,d\bfx dt 
\,+\,\int_{0}^{T}\int_{\Omega}\left[\bfq_E \bfE_a\,+\,\bfq_H\bfH_a
\right]\,d\bfx dt 
\end{displaymath}
\begin{displaymath}
\fl
=\, \int_{0}^{T}\int_{\Omega}\bigg[ 
-\bfE_{f}\epsilon \frac{\partial\bfE_{a}}{\partial t}
\,-\, \bfH_{f} \ncurl \bfE_{a}
\,+\, \bfE_f \sigma \bfE_a
\,-\,\bfH_{f} \mu \frac{\partial \bfH_{a}}{\partial t}
\,+\, \bfE_{f} \ncurl \bfH_{a}\bigg]\,d\bfx dt 
\end{displaymath}
\begin{displaymath}
\fl
+\,\int_{0}^{T}\int_{\Omega}\left[\epsilon \bfE_{f} \phi\,+\,
\mu \bfH_f \psi\right]\,d\bfx dt 
\,+\,\int_{0}^{T}\int_{\Omega}\left[\bfq_E \bfE_a\,+\,\bfq_H\bfH_a
\right]\,d\bfx dt 
\end{displaymath}
\begin{equation}
 \label{Pt4.1}
\fl
+\,\int_{0}^{T}\int_{\partial\Omega}\bfE_a\times\bfH_f\cdot \bfn\,d\sigma dt\,+\,
\int_{0}^{T}\int_{\partial\Omega}\bfE_f\times\bfH_a\cdot \bfn\,d\sigma dt
\end{equation}
\begin{displaymath}
\fl
\,+\,\int_{\Omega}\epsilon\big[(\bfE_{f} \bfE_{a})(T)-
(\bfE_{f} \bfE_{a})(0)\big]\,d\bfx 
\,+\,\int_{\Omega}\mu\big[ (\bfH_{f}\bfH_{a})(T)-
(\bfH_{f}\bfH_{a})(0)\big]\,d\bfx
\end{displaymath}
This equation has the form (\ref{Pt2.1}).
Notice that we have augmented Green's formula
in  (\ref{Pt4.1}), as already shown in (\ref{Pt2.1}),
  by some terms which appear in identical
form  on the left hand side and on the right hand side.

The first term on the left hand side of  equation (\ref{Pt4.1}) 
and the third term on the right hand side cancel each other
 because of (\ref{Etrm.1}) and (\ref{Etrm.2}). The second term
on the left hand side and the first term on the right hand side cancel each
other because of (\ref{Etrm.4}), (\ref{Etrm.5}). The ($t=0$)-terms 
 and the  ($t=T$)-terms vanish due to (\ref{Etrm.3}) and (\ref{Etrm.6}).
 The boundary terms vanish because
of zero boundary conditions. We are left over with the equation
\begin{displaymath}
\int_{0}^{T}\int_{\Omega}\left[ \epsilon\bfE_{f} \phi\,+\,
\mu\bfH_f \psi\right]\,d\bfx dt 
\,=\,
\int_{0}^{T}\int_{\Omega}\left[\bfq_E \bfE_a\,+\,\bfq_H\bfH_a
\right]\,d\bfx dt .
\end{displaymath}
Defining $F^{*}$ by (\ref{Etrm.7}), this can be written as
\begin{displaymath}
\left\langle F \left( {\bfq_E \atop \bfq_H}\right)\,,\, 
\left( {\phi \atop \psi}\right)\right\rangle_U \,=\,
\left\langle  \left( {\bfq_E \atop \bfq_H}\right)\,,\, 
F^{*}\left( {\phi \atop \psi}\right)\right\rangle_U. 
\end{displaymath}
Therefore, $F^{*}$ is in fact the adjoint of $F$, and the lemma 
is proven. \hfill $\Box$\\

\clearpage 
 
\end{document}